  \documentclass[authoryear,preprint,3p,singlecolumn,11pt]{elsarticle}

\usepackage{xcolor}
\usepackage{amsthm}
\usepackage{subfigure}
\usepackage{graphicx}
\usepackage{caption}
\usepackage{amssymb}
\usepackage{lineno}
\usepackage{epsfig}
\usepackage[fleqn]{amsmath}
\usepackage{verbatim}
\newcommand{\eat}[1]{}

\newtheorem{Lemma}{Lemma}

\usepackage{tabulary}
\usepackage{booktabs}
\usepackage[ruled]{algorithm2e}
\usepackage{adjustbox}
\usepackage{rotating}
\usepackage{multirow}
\usepackage{lscape}
 \usepackage{setspace}
\usepackage{url}
\usepackage{xcolor}
\usepackage{soul}

\usepackage[colorlinks, linkcolor=blue,citecolor=blue]{hyperref}

\DeclareMathOperator{\R}{\mathbb{R}}
\makeatletter
\def\BState{\State\hskip-\ALG@thistlm}
\makeatother

\makeatletter
\def\ps@pprintTitle{%
   \let\@oddhead\@empty
   \let\@evenhead\@empty
   \def\@oddfoot{\reset@font\hfil\thepage\hfil}
   \let\@evenfoot\@oddfoot
}
\makeatother

\begin{document}
\begin{frontmatter}

\title{Decomposition approach for Stackelberg P-median problem with user preferences}

 \author[NTU]{Qingyun~Tian}

  \author[NUS]{Yun~Hui~Lin\corref{cor1}}
 \ead{isemlyh@gmail.com}

  \author[NUS]{Dongdong He}

 \address[NTU]{School of Civil and Environmental Engineering, Nanyang Technological  University}
 \address[NUS]{Department of Industrial Systems Engineering and Management, National University of Singapore}
 \cortext[cor1]{Corresponding author.}

\begin{abstract}

The P-median facility location problem with user preferences (PUP) studies an operator that locates P facilities to serve customers/users in a cost-efficient manner, upon anticipating customer preferences and choices. The problem can be visualized as a leader-follower game in which the operator is the leader that opens facilities, whereas the customer is the follower who observes the operator’s location decision at first and then seeks services from the most preferred facility. Such a modeling perspective is of practical importance as we have witnessed its applications to various problems, such as the establishment of power plants in energy markets and the location of healthcare service centers for COVID-19 Vaccination. Despite that a considerable number of solution methodologies have been proposed, many of them are heuristic methods whose solution quality cannot be easily verified. Moreover, due to the hardness of the problems, existing exact approaches have limited performance. Motivated by these observations, we aim to develop an efficient exact algorithm for solving large-scale PUP models. We first propose a branch-and-cut decomposition algorithm and then design accelerated techniques to further enhance the performance.  Using a broad testbed, we show that our algorithm outperforms various exact approaches by a large margin, and the advantage can go up to several orders of magnitude in terms of computational time in some datasets. Finally, we conduct sensitivity analysis to draw additional implications and to highlight the importance of considering user preferences when they exist.

\end{abstract}

\end{frontmatter}

\section{Introduction}
\label{intro}

The classical facility location problem (FLP) studies the strategic planning of an operator who provides services to customers. The problem involves the determination of a set of open facilities from candidate sites, aiming to minimize the cost of satisfying customer's demand/or maximize some objectives of interest (e.g., to maximize the profit of the operator or to cover as many customers as possible).  To date, a large number of models have been proposed to address the FLPs under various settings. For comprehensive discussions of models, algorithms, and applications, we refer readers to~\citet{melo2009facility} and \citet{daskin2011network}.

In the FLP literature, it is commonly assumed that customers (users) are allocated to their closest facilities or the operator acts as a \textit{centralized} decision-maker that assigns customers to facilities and serves them in a cost-efficient fashion. However, this allocation/assignment rule is relatively simple and does not account for customer preferences.  In many scenarios, customers may not prefer the assigned facility. They may instead seek services from some other facility, from which it could be costly for the operator to serve them. This implies that there could be a mismatch between the preferences of the customers and the desirable customer allocation/assignment of the operator. Therefore, when locating the facilities, the operator needs to consider this mismatch, giving rise to a variant of the FLP that is often referred to as the facility location problem with user preferences (FLPUP) or with order~\citep{camacho2014solving,hansen2004lower,maric2012metaheuristic,vasil2009new,vasilyev2010branch}. When the number of facilities to be set up is fixed, the problem is called the P-median problem with user preferences (PUP)~\citep{casas2017solving,camacho2014p,diaz2017grasp}.

\begin{figure*}[h]
\begin{center}
     \psfig{figure=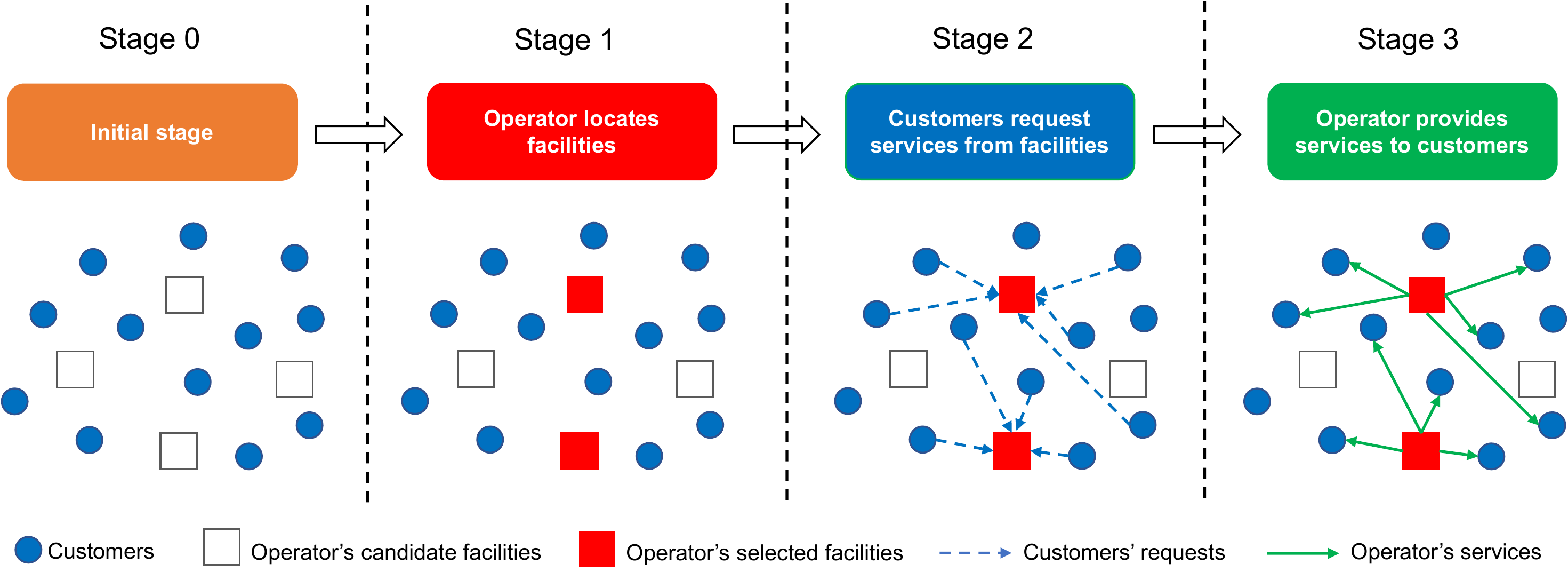,width=165mm}
     \caption{Decision stages of the facility location problem with user preferences.}\label{fig:decision}
\end{center}
\end{figure*}

In the  context of FLPUP and PUP, customers are studied as individual decision-makers that select facilities to serve them based on their preferences. The resultant problem is visualized from a fully \textit{decentralized} perspective and thus can be characterized as a Stackelberg leader-follower game, in which the operator is the leader and the customer is the follower. Figure~\ref{fig:decision} illustrates the decision process of both the operator and the customer. At the initial stage, the customer and the candidate facilities are known as discrete nodes. At Stage 1, the operator selects locations from the candidate sites. Customers, as followers, observe the operator's decisions at first and then request services from facilities at Stage 2. Finally, the operator proceeds to serve customers according to the requests. 

It is worth noting that when locating facilities at Stage 1, the operator needs to anticipate how customers choose facilities at Stage 2 because the decentralized customers want to be served by the most preferred facility that is not necessarily cost-efficient from the operator's side. Customers' requests could thus have a significant impact on the operators' service profile and the service provision cost at Stage 3. Since facilities are generally built for a long-term purpose and it will be costly and difficult to alter the decision, it is important for the operator to consider customer preferences and anticipate customers' choices. Such a modeling perspective is of practical importance. Recently, we have witnessed several applications of the FLPUP and PUP to real-world problems, such as the establishment of power plants in energy markets \citep{lotfi2021robust} and the location of healthcare service centers for COVID-19 Vaccination~\citep{cabezas2021two}. 

This paper investigates a PUP model.  Due to the hidden bilevel structure,  the resultant mathematical model is computationally challenging. Until now, we are unaware of any existing exact solution approach that can solve large-scale PUP instances within a reasonable computational time. This paper thus aims to develop an efficient exact algorithm for large-scale PUP instances. Below, we briefly review the related literature with a focus on the solution approaches to clearly position our contributions.  \\

\noindent \textbf{Literature review.} The first study that considered customer preferences is~\cite{hanjoul1987facility} where each customer has a preference ordering towards the set of potential locations, and the simple plant location problem was presented with a set of constraints on preference ordering. A greedy heuristic algorithm was proposed to conduct the numerical tests on multiple instances.   Since then, there is a proliferation of research on the FLPUP and PUP.  

Although the problem has a Stackelbeg structure (as illustrated in Figure~\ref{fig:decision}), single-level formulation indeed exists. A typical exact approach thus consists of formulating the FLPUP and PUP as a single-level mixed-integer linear programming (MILP) and then solving the resultant MILP using existing solvers (possibly strengthened by some valid inequalities). In general, the single-level formulation can be divided into two research streams. The first one is based on the primal-dual information of the lower-level problem and leverages the KKT conditions to recast the bilevel model into a mixed-integer program with complementarity constraints. Then the non-convex complementarity constraints are linearized, giving rise to a MILP that can be directly solved by off-the-shelf solvers~(e.g., see~\cite{camacho2014p,casas2017solving,cao2006capacitated,lotfi2021robust}). However, due to the use of the KKT conditions, such an approach needs to introduce a large number of variables (i.e., the dual variables) that dramatically increases  the formulation size. As a consequence, its performance is typically limited. The second stream, essentially, utilizes the \textit{closest assignment constraints} (CACs) to subtly rewrite the lower-level problem into a set of linear constraints without additional variables. For a comprehensive discussion on the CACs, we refer readers to \cite{espejo2012closest}. To date, a considerable number of studies have adopted the CAC-based reformulation. For example, \cite{canovas2007strengthened}  proposed a strengthened CAC to  speed up solving the model of~\cite{hanjoul1987facility}.  \cite{vasil2009new} provided several valid inequalities when the bilevel problem was reformulated into a MILP using CACs and demonstrated the tightness of the inequalities in terms of their linear relaxations and integrality gaps. Then \cite{vasilyev2010branch} implemented a branch-and-cut method based on the proposed family of valid inequalities in \cite{vasil2009new}. A recent application can be found in~\cite{casas2017solving}.

We point out that the above approaches are exact and can yield the proven optimal solution, provided that the reformulated MILP model can be solved optimally. However, the MILP model itself is a challenging problem that cannot scale well on the problem size, rendering these approaches computational prohibitive for large-scale instances. Therefore, heuristic approaches are gaining popularity.  Recent years have witnessed the use of Lagrangian-based heuristics to various problems involving customer preferences~(e.g., see \cite{cabezas2018lagrangean,cabezas2021two,lee2012facility}). Moreover, there is an explosion of the related research on metaheuristics.  \cite{maric2012metaheuristic} proposed three metaheuristic methods for solving FLPUP, i.e., Particle Swarm Optimization (PSO), Simulated Annealing (SA), and a combination of Reduced and Basic Variable Neighborhood Search Method (RVNS-VNS). To solve a large-scale bilevel FLPUP, \cite{camacho2014solving} developed a population-based evolutionary algorithm (EA). A similar EA was applied to solve a bilevel facility location problem with cardinality constraints and preferences~\citep{calvete2020matheuristic}. \cite{casas2018approximating} obtained approximated solutions to the problem by using a cross-entropy method to search for improved location decisions. In \cite{casas2017solving}, the authors investigated the PUP model and designed a scatter search metaheuristic algorithm to efficiently handle large-scale problems. In the context of maximum coverage, \cite{diaz2017grasp} presented  a bilevel model considering that customers will choose the facilities within a coverage radius. They proposed a greedy randomized adaptive search procedure (GRASP) heuristic and a hybrid GRASP-Tabu heuristic to find near-optimal solutions. Recently, \cite{mrkela2021variable} investigated a similar problem as  in \cite{diaz2017grasp}. Instead of locating P facilities, they considered a limited budget for locating facilities, and an efficient Variable Neighborhood Search (VNS) was proposed as the solution approach.  \\

\noindent \textbf{Our contributions.}  Despite that there is a considerable number of research works dedicated to developing solution methodologies for the FLPUP and PUP, many of them are heuristic methods whose solution quality cannot be easily verified. Moreover, due to the hardness of the problems,  existing exact approaches based on single-level reformulation techniques have limited performance. Such an observation motivates this paper to develop an efficient exact algorithm that can be applied to solve large-scale PUP models. In summary, our contributions are three-fold.

 (i)\textit{A branch-and-cut Benders decomposition algorithm}. By exploring the problem structure, we prove that in a CAC-based MILP model~\citep{casas2017solving}, the high-dimensional binary variables used to model customer preferences and facility choices can be relaxed to continuous variables. Using this condition, we propose a branch-and-cut algorithm, where we project out the preference-related variables and solve a master problem defined on a low-dimensional decision space. 

(ii) \textit{An enhanced Benders approach with analytical separation}. The decomposition algorithm is expected to be more efficient than the MILP model since the branch-and-cut searching tree only involves the location variable of the original MILP and the problem size is thus substantially smaller.  However, our numerical experiment reveals that such an algorithm cannot effectively expedite the solution process because generating the Benders cut (i.e., Benders separation) at nodes of the tree is not trivial, and significant computational time can be spent on the separation function. In light of this, we design an analytical approach for Benders separation to further speed up the algorithm. It turns out that our specialized separation is so efficient that we observe a significant performance improvement (with up to orders of magnitude in terms of computational time). Using a broad testbed, our computational experiments support that the enhanced Benders approach is able to handle large-scale instances satisfactorily. 

(iii) \textit{Additional managerial implications}. We  conduct sensitivity analysis using instances from a standard benchmark library. Our analysis reveals  that when user preferences indeed exist (such as in the context of energy markets \citep{lotfi2021robust} and healthcare service~\citep{cabezas2021two}, the operator must consider these preferences and correctly anticipate user’s choices; otherwise, the operator could suffer from a substantially higher cost, and opening more facilities could ironically result in additional service costs. \\

The rest of the paper is organized as follows. We introduce the problem and the model in Section~\ref{S:pd}. We then develop the branch-and-cut Benders decomposition and the acceleration technique in Section~\ref{S:sm}. Extensive numerical studies are conducted in Section~\ref{S:ns} to demonstrate the effectiveness of our algorithm and provide additional managerial implications. Finally, we conclude the paper in Section~\ref{S:cl}.

\section{Model formulation}
\label{S:pd}

This section describes the model formulation of the P-median problem with user preferences (PUP).  Table~\ref{tab:Nomenclature} provides the main notations.
\begin{table}[h]
\caption{Notations}
\label{tab:Nomenclature}       
\begin{tabular}{lll}
\hline\noalign{\smallskip}
\multicolumn{3}{l}{\underline{Sets}}\\
$I$ &      $:$ & set of customers.\\
$J$ &     $:$ & set of candidate facilities.\\

\multicolumn{3}{l}{\underline{Parameters}}\\
$c_{ij}$ &     $:$ & cost of serving customer $i$ by facility $j$, $\forall i \in I, j \in  J$\\
$g_{ij}$ &     $:$ & disutility of facility $j$ to customer $i$, $\forall i \in I, j \in J$\\
$\pi_{ij}$ &     $:$ & normalized disutility of facility $j$ to customer $i$, defined as $g_{ij}/(\max_{j\in J} g_{ij})$,  $\forall i \in I, j \in J$\\
$P$ & $:$ & the number of open facilities \\

\multicolumn{3}{l}{\underline{Decision Variables}}\\
$x_j$ & $:$ &  binary. If facility $j$ is open; 0, otherwise, $\forall j \in J$ \\
$y_{ij}$ & $:$ & binary. If customer $i$ seeks the service from facility $j$; 0, otherwise, $\forall i \in I, j \in J$ \\
\noalign{\smallskip}\hline
\end{tabular}
\end{table}

\subsection{Stackelberg visualization}

The virtual decision process of [PUP] is as explained in Figure~\ref{fig:decision}. We use $J$ to denote the set of candidate facilities and $I$ to denote the set of customers. At Stage 1, the operator will open $P$ facilities, selecting from set $J$. We define a binary variable $x_j, \forall j \in J$, which is 1 if facility $j$ is open; 0, otherwise. After $P$ facilities are built, customers then determine which facilities to seek services from at Stage 2. To reflect customers' choices, we define a binary variable $y_{ij}$ such that $y_{ij} = 1$ if customer $i$ requests services from facility $j$, $\forall i \in I, j \in J$. 

Following the standard assumption in the FLPUP and the PUP literature, we assume that facility $j$ has a disutility to customer $i$, i.e., $g_{ij}$, and customer $i$ prefers facility $j$ over facility $k$, if $g_{ij} < g_{ik}$. In other words, when both facility $j$ and facility $k$ are open, customer $i$ will not select facility $k$ if $g_{ij} < g_{ik}$. With these definitions, given a location decision $x$, customer preferences to the facilities can be modeled by
\begin{subequations}\label{model:stage2}
\begin{alignat}{10}
\label{stage2:1}y^*= \arg\min_{y}~&\sum_{i \in I}\sum_{j \in J} g_{ij}y_{ij} \\
\text{st.}~&\sum_{j \in J} y_{ij} = 1 &\forall i \in I\\
\label{stage2:2} &y_{ij} \leq x_j  \quad &\forall i \in I, j \in J \\
\label{stage2:3}& y_{ij} \in \{0,1\} \quad &\forall i \in I, j \in J
\end{alignat}
\end{subequations}
which describes the decision problem at Stage 2. The objective function imposes that the customer will minimize the disutility. This is equivalent to state that he/she will request the service from the most preferred facility located by the operator. 

Given the solution $y^*$ from the above problem, the operator incurs a service provision cost of $\phi = \sum_{i \in I}\sum_{j\in J} c_{ij}y^*_{ij}$ at Stage 3, where $c_{ij}$ is the cost of serving customer $i$ from facility $j$ (possibly weighted by the demand size or the relative ``importance" of the customer). Note that the operator also wants to minimize the total service cost; therefore, when locating facilities at Stage 1, the operator needs to solve the following problem
\begin{subequations}\label{model:P}
\begin{alignat}{10}
\label{P:cons1} \min~ \phi = & \sum_{i \in I} \sum_{j \in J} c_{ij}y^*_{ij} \\
\textbf{[PUP]} \quad \qquad \label{P:cons2}\text{st.}~&\sum_{j \in J} x_{j} = P \\
\label{P:cons3}& x_{j} \in\{0,1\} \quad \forall j \in J\\
\label{P:cons4}& y^* \in (\ref{model:stage2}) 
\end{alignat}
\end{subequations}
where Constraint~(\ref{P:cons4}) states that $y^*$ is obtained by solving Problem~(\ref{model:stage2}). Therefore, [PUP] is indeed a \textit{mixed-integer bilevel program}. \\

\noindent \textbf{Remark 1.}  Through this paper, we assume that customer preferences for facilities are different, i.e., elements in the vector $g_{i\cdot}$ are distinct ($g_{ij} \ne g_{ik}$ if $j \ne k$). Under this assumption, an operator's decision will only lead to a unique solution of $y$ (i.e., Problem~\ref{model:stage2} has an unique solution), allowing us to avoid the discussion on optimistic and pessimistic strategies of the follower in the context of bilevel optimization.

\subsection{Single-level reformulation}

The bilevel structure of [PUP] renders the problem computationally challenging. Here, we present a single-level reformulation model of [PUP] to circumvent the difficulties of handling Problem~(\ref{model:stage2}) in [PUP].   In fact, minimizing the objective function~(\ref{stage2:1}) in Problem~(\ref{model:stage2}) can be represented by  the following inequality
\begin{align}
\label{closest_assign}\sum_{j \in J}g_{ij}y_{ij} \leq g_{ij}x_j + G_i (1-x_j) \quad \forall i \in I, j \in J
\end{align}
where $G_i = \max_{j \in J} g_{ij}$. To see the equivalence, if facility $j$ is not open ($x_j=0$), then the right hand side becomes $G_i$, meaning that  (\ref{closest_assign}) is inactive since $G_i = \max_{j \in J} g_{ij}$. If  facility $j$ is open ($x_j=1$),  then (\ref{closest_assign}) imposes that, among all the possible allocations, customer $i$ will be allocated to the facility with the minimum value of $g_{i\cdot}$. Essentially, this equation belongs to a type of \textit{closest assignment constraints} (CACs). It was proposed by \cite{berman2009optimal} and has been used by~\cite{casas2017solving} in a PUP model.  

Now, define $\pi_{ij} = g_{ij}/G_i$, i.e., the vector $g_{i\cdot}$ is normalized by the large value in it to get $\pi_{i\cdot}$. We have the following single-level reformulation model for [PUP]
\begin{subequations}\label{model:SRM}
\begin{alignat}{10}
\label{SRM:cons1} \min~& \sum_{i \in I} \sum_{j \in J} c_{ij}y_{ij} \\
\label{SRM:cons2}\text{st.}~& \sum_{j \in J} y_{ij} = 1 \quad &\forall i \in I \\
\label{SRM:cons3}\textbf{[SRM]} \qquad \quad  & y_{ij} \leq x_j  \quad &\forall i \in I, j \in J \\
\label{SRM:cons4}& \sum_{j \in J}\pi_{ij}y_{ij} \leq (\pi_{ij}-1)x_j + 1 \quad &\forall i \in I, j \in J \\
\label{SRM:cons6}& y_{ij} \in \{0,1\} \quad &\forall i \in I, j \in J \\
\label{SRM:cons5}& x \in \Omega
\end{alignat}
\end{subequations}
where $\Omega$ is the decision space for the operator
\begin{align}
\Omega = \left\{x \mid \sum_{j \in J} x_j = P, x_j \in \{0,1\}, \forall j \in J \right\}
\end{align}
This set is defined to simplify the notation for our later discussions. Now, [SRM] is a mixed-integer linear program (MILP) and is ready to be solved by modern MILP solvers. \\

\noindent \textbf{Remark 2.} We point out that other types of CACs that can possibly be used to reformulate [PUP] into a single-level model have been discussed in detail in \cite{espejo2012closest}. Among all CACs, we choose (\ref{SRM:cons4}) since it only involves $O(|I|\cdot |J|)$ constraints, whereas others may contain up to $O(|I|\cdot |J| \cdot |J|)$ constraints. \\ 

\noindent \textbf{Remark 3.} Besides the CAC-based reformulation,  the \textit{primal-dual reformulation} has also been discussed in the literature~ \citep{casas2017solving,casas2018approximating,lotfi2021robust}. The fundamental idea is to replace Problem~(\ref{model:stage2}) with its KKT conditions, thereby recasting [PUP] into a single-level primal-dual reformulation model [PDRM]. We provide the  [PDRM] formulationl in~\ref{app:pdr}. We note that this approach introduces a large number of new variables into the formulation (i.e., there are $|I|\cdot(2|J|+1)$ dual variables in the KKT conditions)  and thus significantly increases the problem size. Nevertheless, the single-level reformulation approach leveraging the KKT conditions of the lower-level problem is probably the most widely used approach for solving bilevel programs; therefore, we will keep [PDRM] in this paper and use it as one of the benchmarks to our proposed algorithms.

\section{Solution methodology}
\label{S:sm}

This section presents the solution methodology for [SRM]. In the literature, the single-level MILP model is typically solved using off-the-shelf solvers (e.g., CPLEX and Gurobi). As shown in~\cite{casas2017solving},  such a straightforward approach cannot handle large-scale problems. Therefore, we develop a branch-and-cut Benders decomposition algorithm to expedite solving [SRM], where Benders separation (i.e., the procedure of generating Benders cuts) at integer nodes of the searching tree is performed leveraging external solvers. However, we observe that such a standard Benders approach is still not efficient enough. Therefore, we further derive an acceleration technique to enhance the algorithm performance.

\subsection{Standard branch-and-cut Benders decomposition}
\label{SS:standard}

Our decomposition approaches are inspired by the following observation.
\begin{Lemma} \label{lemma_relax}
In the  [SRM] formulation, Constraint~(\ref{SRM:cons6}) can be relaxed to $y_{ij} \geq 0$ without changing the solution and the objective.
\end{Lemma}
\begin{proof}
Suppose a location decision is made at $\bar{x}$. Let $\bar{\tau}$ be the set of open facilities, i.e., $\bar{\tau} = \left\{ \forall j \in J \mid \bar{x}_j = 1 \right\}$. Based on the condition $y_{ij} \leq \bar{x}_j$, Constraint~(\ref{SRM:cons4}) can be restated as $ \sum_{j \in \bar{\tau} }\pi_{ij}y_{ij} \leq \min_{j \in \bar{\tau}} \pi_{ij}, \forall i \in I$. For notation simplicity, we drop subscript $i$ in this proof. 

Let $m$ denote the most preferred open facility for customer, i.e., $m =  \arg \min_{j\in \bar{\tau}} \pi_{j}$. We have $ \sum_{j \in \bar{\tau} }\pi_{j}y_{j} \leq \pi_{m}$. Note that $\pi_m < \pi_j, \forall j \in  \bar{\tau} \setminus \{m\}$, and $\sum_{j \in \bar{\tau}} y_j = 1$. Therefore, $ \sum_{j \in \bar{\tau} }\pi_{j}y_{j} \leq \pi_{m}$ holds only if $y_m = 1$ and $y_j = 0, \forall j \in \bar{\tau} \setminus \{m\}$, which also uses the condition that the elements in $\pi_{i\cdot}$ ($g_{i\cdot}$) are distinct.  This indicates that relaxing $y_j \in \{0,1\}$ to $ y_j \geq 0$ will not change the solution since $y_j$ in [SRM] will automatically be an integer.
\end{proof}
This lemma implies that when the operator's decision is made, the remaining subproblem of [SRM] will become a linear program (LP) that prosesses strong duality. Therefore, it is possible to design a decomposition algorithm leveraging the dual information of the subproblem, which is exactly the idea of Benders decomposition. Modern Benders decomposition algorithms are typically implemented within the solver's  branch-and-cut framework, and we have observed a considerable number of research works successfully employing such a technique to address facility location problems~\citep{cordeau2019benders,fischetti2016benders,ljubic2012exact,taherkhani2020benders}.  

In this paper, we design a branch-and-cut Benders decomposition to solve [SRM]. Our approach only maintains the location decision $x$ in a master problem by projecting out the continuous variable $y$ that is used to reflect the preferences and compute the service cost. 
Specifically, the master problem is defined as
\begin{subequations}\label{model:MP}
\begin{alignat}{10}
\label{MP:obj} \min~& \sum_{i \in I} w_i \\
\label{MP:cons1}\textbf{[MP]} \qquad \text{st.}~& w_i \geq \Phi_i(x) \quad \forall i \in I \\
\label{MP:cons2}& x \in \Omega
\end{alignat}
\end{subequations}
where $\Phi_i(x)$ is obtained by solving the following subproblem
\begin{subequations}\label{model:SP}
\begin{alignat}{10}
\label{SP:obj} \Phi_i(x) =  \min~& \sum_{j \in J} c_{ij}y_{ij} \\
\label{SP:cons1}\text{st.}~ & \sum_{j \in J} y_{ij} = 1  & (\mu_i) \\
\label{SP:cons2}\textbf{[SP$_i$$(x)$]} \qquad \qquad \qquad & y_{ij} \leq x_j  \quad \forall j \in J & (\lambda_{ij})  \\
\label{SP:cons3}& \sum_{j \in J}\pi_{ij}y_{ij} \leq (\pi_{ij}-1)x_j + 1 \quad \forall j \in J~~& (v_{ij})  \\
\label{SP:cons4}& y_{ij} \geq 0 \quad \forall j \in J
\end{alignat}
\end{subequations}
Note that $\Phi_i(x)$ is a convex function over $x$; therefore, [MP] is a mixed-integer convex optimization problem, and we can solve it by approximating $\Phi_i(x)$ with its linear under-estimators, i.e., the Benders cuts.

Given an integer solution $\bar{x}$ from [MP], the Benders cut can be obtained by solving the dual problem of [SP$_i$$(\bar{x})$]. Let $\mu_i$ and $\lambda_{ij}$ and $v_{ij}$ in (\ref{model:SP}) be the dual variables associated with the constraints. The dual subproblem is given by
\begin{subequations}\label{model:DSP}
\begin{alignat}{10}
\label{DSP:obj} \max~& \mu_i - \sum_{j\in J}v_{ij} + \sum_{j \in J}(v_{ij} - \pi_{ij} v_{ij} - \lambda_{ij}) \bar{x}_j \\
\label{DSP:cons1}\textbf{[DSP$_i$$(\bar{x})$]} \qquad \qquad \text{st.}~ & c_{ij} - \mu_i + \lambda_{ij} + \pi_{ij} \sum_{j \in J} v_{ij} \geq 0 \quad \forall  j \in J \\
\label{DSP:cons2}& \lambda_{ij} \geq 0\quad \forall  j \in J \\
\label{DSP:cons3}&v_{ij} \geq 0\quad \forall  j \in J
\end{alignat}
\end{subequations}
Let $(\bar{\lambda},\bar{v},\bar{\mu})$ be the optimal solution to the above LP. The following Benders cut defined at $\bar{x}$ arises.
\begin{align}\label{benders_cuts}
 w_i \geq \Phi_i(x) \geq \bar{\mu}_i - \sum_{j\in J}\bar{v}_{ij} + \sum_{j \in J}(\bar{v}_{ij} - \pi_{ij} \bar{v}_{ij} - \bar{\lambda}_{ij}) x \quad \forall i \in I
 \end{align}
We can now solve [MP] using the following MILP formulation, i.e., the \textit{relaxed master problem},
\begin{subequations}\label{model:MP}
\begin{alignat}{10}
\label{MP:obj} \min~& \sum_{i \in I} w_i \\
\label{MP:cons1}\textbf{[rMP]} \qquad \text{st.}~& w_i \geq \bar{\mu}_i - \sum_{j\in J}\bar{v}_{ij} + \sum_{j \in J}(\bar{v}_{ij} - \pi_{ij} \bar{v}_{ij} - \bar{\lambda}_{ij}) x \quad \forall i \in I, (\bar{\lambda},\bar{v},\bar{\mu}) \in \Xi \\
\label{MP:cons2}& x \in \Omega
\end{alignat}
\end{subequations}
where $\Xi $ is the set of dual variables. We can expand the set, each time we have an integer solution $\bar{x}$ and solve the corresponding [DSP$_i$$(\bar{x})$]. These dual variables define Benders cuts that cut off non-optimal solutions and lead us to the proven optimal solution. 

We then illustrate our implementation of Benders decomposition. Modern advanced solvers, such as CPLEX and Gurobi, use branch-and-cut algorithms to solve (mixed-)integer programs. Apart from the general-purpose cuts that are embedded within the solvers, we can design a problem-specific separation function that uses a  [rMP] solution $\bar{x}$ as an input and generates violated cuts, which are inserted into the searching tree to improve the relaxation bound and/or cut off non-optimal solutions. We perform Benders separation only at integer nodes because this is sufficient to guarantee the optimal convergence. More specifically, we use Gurobi 9.1.2 and rely on its default branching rules and embedded cutting planes. We initialize the master problem [rMP] without Benders cuts  (i.e., $\Xi = \emptyset$). Now, when Gurobi visits an integer node, we retrieve the integer value of $\bar{x}$ and generate the corresponding Benders cuts as in (\ref{benders_cuts}) by solving [DSP$_i$$(\bar{x})$].  The Benders separation is performed within the \textit{callback} function of Gurobi. Cuts are added as \textit{lazy cuts} to the current node relaxation only if they violate the current solution $(\bar{x},\bar{w})$ by more than $10^{-5}$ of the relative violation, defined as the absolute violation of the cut divided by the current $\bar{w}_i$ value. Due to the convexity of $\Phi_i$, these inserted cuts are globally valid and will eventually lead us to the optimal solution of [SRM] with a zero MIP gap.

\subsection{Accelerating Benders approach through analytical separation}
\label{SS:accelerating_Benders}

The above branch-and-cut Benders decomposition generates Benders cuts by solving the dual subproblem [DSP$_i(\bar{x})$] using external solvers, each time an integer node is visited during the searching tree. However, [DSP$_i(\bar{x})$] itself is a large-scale LP. For some problems, solving [DSP$_i(\bar{x})$] is not trivial and may take up to several seconds. Noting that we typically need to visit a large number of integer nodes before the MIP gap vanishes, significant computational time could thus be consumed on Benders separation, causing a bottleneck in the branch-and-cut algorithm. In light of this, we design a specialized analytical approach for Benders separation, which is very fast and can, in general, expedite the algorithm by a large margin as demonstrated in our experiment. One should keep in mind that the procedure explained below relies on the condition that the solution $\bar{x}$ is binary. \\

\noindent \textbf{Step 1: Solve [SP$_i$$(\bar{x})$]}. To start with, we point out that given a location decision $\bar{x}$, [SP$_i$$(\bar{x})$] can be solved using a simple sorting algorithm since customers will seek the service from the most preferred open facility. Let $\bar{\tau}$ denote the set of open facilities, i.e.,
\begin{align}
\bar{\tau} = \left\{ \forall j \in J \mid \bar{x}_j = 1 \right\}
\end{align}
Let $m$ denote the most preferred open facility for customer $i$. We have
\begin{align}
m =  \arg \min_{j\in \bar{\tau}} \pi_{ij}
\end{align}
In other words, given $\bar{x}$, customer $i$ will select facility $m$, and thus the optimal solution of $y$ is
\begin{equation}
\bar{y}_{ij} =
 \begin{cases}
1, \text{if}~j = m \\
0, \text{otherwise}
\end{cases}
\end{equation}
Then, the leader incurs a cost of $\sum_{j \in J}c_{ij}\bar{y}_{ij} = c_{im}$ to serve customer $i$, meaning that the optimal value of [SP$_i$$(\bar{x})$] is
\begin{align}
\label{eqt:value}\Phi_i(\bar{x}) = c_{im}
\end{align}
Upon here, we have solved [SP$_i$$(\bar{x})$] and obtained the objective $c_{im}$ and the solution $\bar{y}$. \\

\noindent \textbf{Step 2: Transform [DSP$_i$$(\bar{x})$]}. To solve [DSP$_i$$(\bar{x})$], we first transform it into a problem of finding feasible solutions to a set of linear equations. Based on the value of $\bar{x}_j$ and $\bar{y}_{ij}$, there exist three cases.
\begin{description}
\item Case 1: $\bar{x}_j = 1$ and $\bar{y}_{ij} = 0$.  In [SP$_i$$(\bar{x})$], constraint~(\ref{SP:cons2}) is inactive; therefore, $\lambda_{ij} = 0$. Moreover, constraint~(\ref{SP:cons3}) becomes $\pi_{im} \leq \pi_{ij}$, which holds since $m =  \arg \min_{j\in \bar{\tau}} \pi_{ij}$. This indicates that constraint~(\ref{SP:cons3}) is also inactive, leading to $v_{ij} = 0$.

\item Case 2: $x_j = 0$ and $y_{ij} = 0$. We have $v_{ij} = 0$ since constraint~(\ref{SP:cons3}) is inactive.

\item Case 3: $x_j = 1$ and $y_{ij} = 1$. We have $j = m$.
\end{description}
Combining these cases, we conclude that $\sum_{j \in J}v_{ij} = v_{im}$ since if $j \neq m$, then $v_{ij} = 0$. Therefore, the objective function under $\bar{x}$ becomes
\begin{align}
\mu_i - v_{im} - \sum_{j \in \bar{\tau} }(v_{ij} - \pi_{ij} v_{ij} - \lambda_{ij}) =  \mu_i - \pi_{im} v_{im}  - \sum_{j \in \bar{\tau} }\lambda_{ij} = \mu_i  - \lambda_{im} -  \pi_{im} v_{im}
\end{align}
The last equality holds since when $j \in \bar{\tau}$, $\lambda_{ij}$ can be non-zero only if $j = m$ according to Case 1.

Note that, for customer $i$, Case 1 reflects the dual information of the set of open facilities excluding the most referred one, i.e.,  $\bar{\tau} \setminus \{m\}$, whereas Case 2  reflects the dual information of the set of facilities are not open, i.e, $J \setminus \bar{\tau}$. We can therefore rewritten [DSP$_i$$(\bar{x})$] as
\begin{subequations}\label{model:DSP-new}
\begin{alignat}{10}
\max~& \mu_i  - \lambda_{im} -  \pi_{im} v_{im}\\
\text{st.}~ & c_{ij} - \mu_i + \pi_{ij}  v_{im} \geq 0 \quad & \forall j \in \bar{\tau} \setminus \{m\} \\
& c_{ij} - \mu_i + \lambda_{ij} + \pi_{ij}  v_{im} \geq 0 \quad & \forall  j \in J \setminus \bar{\tau} \\
& (\lambda,v) \geq 0
\end{alignat}
\end{subequations}
Note that the optimal value is $c_{im}$ according to (\ref{eqt:value}). We have  $c_{im} = \mu_i  - \lambda_{im} -  \pi_{im} v_{im}$ based on the strong duality of LP. Therefore, solving Problem (\ref{model:DSP-new}) is equivalent to find out a feasible solution of the following linear equations
\begin{subequations}\label{model:linear_equation}
\begin{alignat}{10}
\label{LE:1} & c_{ij}  + \pi_{ij}  v_{im} \geq  c_{im} + \lambda_{im} + \pi_{im} v_{im} \quad & \forall j \in \bar{\tau} \setminus \{m\} \\
\label{LE:2} & c_{ij}  + \lambda_{ij} + \pi_{ij}  v_{im} \geq  c_{im} + \lambda_{im} + \pi_{im} v_{im}  \quad & \forall  j \in J \setminus \bar{\tau} \\
\label{LE:3} & (\lambda,v) \geq 0
\end{alignat}
\end{subequations}
where we replace $\mu_i$ with $c_{im} + \lambda_{im} +  \pi_{im} v_{im}$. \\

\noindent \textbf{Step 3: Obtain the dual variables}. Now, we proceed to deriving an analytical feasible solution to Problem (\ref{model:linear_equation}). Note that we can safely set $\lambda_{im} = 0$ since to satisfy equations (\ref{LE:1}) and  (\ref{LE:2}), $\lambda_{im}$ should be small enough. After some algebra, (\ref{LE:1}) becomes
\begin{equation}
\begin{aligned}
& v_{im} \geq \frac{c_{ij} - c_{im}}{\pi_{im}-\pi_{ij}} \quad \forall j \in \bar{\tau} \setminus \{m\}
\end{aligned}
\end{equation}
Since $v_{im} \geq 0$, we set
\begin{align}
v_{im} = \max_{j \in \bar{\tau} \setminus \{m\} } \left[ \frac{c_{ij} - c_{im}}{\pi_{im}-\pi_{ij}} \right]_+
\end{align}
where $[z]_+ = \max\{z,0\}$, and compute $\mu_i$ by
\begin{align}
\mu_i =  c_{im} + \pi_{im} v_{im}
\end{align}
which utilizes the condition $c_{im} = \mu_i  - \lambda_{im} -  \pi_{im} v_{im}$ and $\lambda_{im}$ is 0. Furthermore, based on (\ref{LE:2}), we have
 \begin{align}
& \lambda_{ij}  \geq \mu_i - c_{ij} - \pi_{ij}  v_{im} \quad \forall  j \in J \setminus \bar{\tau}
\end{align}
Since $\lambda_{ij} \geq 0$, we set
\begin{align}
\lambda_{ij} =  \left[\mu_i - c_{ij} - \pi_{ij}  v_{im} \right]_+ \quad \forall  j \in J \setminus \bar{\tau}
\end{align}
Note that $\lambda_{ij}=0$ when $\bar{x}_j =1$. We can express $\lambda$  as
\begin{align}
\lambda_{ij} = (1-x_j)\left[\mu_i - c_{ij} - \pi_{ij}  v_{im} \right]_+ \quad \forall j \in J
\end{align}
Finally, we summarize the above three steps compactly into the following lemma.
\begin{Lemma}\label{lemma}
Given an integer solution $\bar{x}$ from [rMP], let $\bar{\tau} = \left\{\forall j \in J \mid \bar{x}_j = 1 \right\}$ and $m =  \arg\min_{j\in \bar{\tau}} \pi_{ij}$. An optimal dual variable ($v,\mu,\lambda$) can be obtained analytically as
\begin{subequations}\label{analytical_formula}
\begin{alignat}{10}
&v_{im} = \max_{j \in \bar{\tau} \setminus \{m\} } \left[ \frac{c_{ij} - c_{im}}{\pi_{im}-\pi_{ij}} \right]_+ \\
&v_{ij} = 0  &\text{if}~j \neq m \\
&\mu_i =  c_{im} + \pi_{im} v_{im} \\
&\lambda_{ij} = (1-x_j)\left[\mu_i - c_{ij} - \pi_{ij}  v_{im} \right]_+ &\forall j \in J
\end{alignat}
\end{subequations}
where $[z]_+ = \max\{z,0\}$.
\end{Lemma}
To conclude the above discussion, given an integer solution $\bar{x}$, the optimal dual variables can be obtained by the procedure described in Lemma~\ref{lemma}, which only involves a sorting algorithm (i.e., determining the most preferred facility) and some elemental matrix manipulation. Therefore, the analytical separation is very efficient. We emphasize that the integrality of $\bar{x}$ is important for the analytical separation because Lemma~\ref{lemma} holds only if $\bar{x}$ is an integer point. Therefore, similar to the standard branch-and-cut Benders decomposition, when we implement the accelerated Benders approach, analytical separation of Benders cuts is performed only at integer nodes.  Moreover, following~\cite{lin2021branch}, we conservatively set the integer feasibility tolerance and primal feasibility tolerance to the minimum values provided by Gurobi, i.e., the Gurobi parameters \textit{IntFeasTol} and \textit{FeasibilityTol} are set to $10^{-9}$ to enhance the numerical stability, despite doing so may potentially increase the overall computational time. Then the current ``integer" $\bar{x}$ from the node is rounded to ensure the integrality.

\section{Numerical study}
\label{S:ns}

This section presents numerical experiments. We first conduct extensive computational experiments using a broad testbed to demonstrate the computational efficiency of the proposed Benders-based approaches. Then we discuss the importance of considering user preferences in facility location problems using  instances from a standard library.

\subsection{Computational experiment}

We start with the computational experiment. Throughout this section, we refer to the standard Branch-and-cut Benders decomposition in Section~\ref{SS:standard} as {\fontfamily{qcr}\selectfont{Benders}} and the accelerated Benders approach with analytical separation in Section~\ref{SS:accelerating_Benders} as  {\fontfamily{qcr}\selectfont{Benders-AS}}. Furthermore, we test the performance of two direct approaches where Gurobi is used to solve [SRM] presented in Section~\ref{S:pd} and  [PDRM] presented in~\ref{app:pdr}. They serve as benchmarks and are referred to as {\fontfamily{qcr}\selectfont{SRM}} and {\fontfamily{qcr}\selectfont{PDRM}}, respectively. All experiments are done on a 16 GB memory macOS computer with a 2.6 GHz Intel Core i7 processor, using Gurobi 9.1.2 as the solver and Python as the programming language.  Moreover, the time limit is set to 7200 seconds.

Our testbed consists of three datasets with different scales and structures.  They are described below.\\

\noindent  \textbf{PMPUP}. A standard testbed for \textit{P-Median Problem with Users Preferences} from the well-known facility location benchmark library \textit{Discrete Location Problems} (see \url{http://www.math.nsc.ru/AP/benchmarks/Bilevel/bilevel-eng.html}).  In total, there are 30 instances (i.e., \textit{inst-333, inst-433, inst-533, ..., inst-3233}). \\ 

\noindent  \textbf{RND}. Medium-scale instances that are randomly generated. Candidate facilities and customers are generated from a two-dimensional uniform distribution $[0,100]^2$. The distance between facility $j$ and customer $i$ (i.e., $l_{ij}$) is represented by the Euclidean distance. Each customer has a demand $d_i$, following a uniform distribution $[0,10]$. The cost of serving customer $i$ from facility $j$ is computed as $c_{ij} = d_i\cdot l_{ij}$. The disutility $g_{ij}$ is randomly generated from  $[(1-\delta)\cdot l_{ij}, (1+\delta)\cdot l_{ij}]$, where $0\leq\delta< 1$. That is, we allow $g_{ij}$ to deviate $100\delta\%$ from $l_{ij}$, and thus, the nearest facility is not necessarily the most attractive facility to customers. Furthermore, we consider the following combinations of problem sizes and parameters: $(|I|,|J|) = \{(200,150),(300,100),(400,150),(500,100)\}$, $P = \{5,10,15,20,30\}$ , $\delta = \{0.3,0.5\}$. We have $5\times 4 \times 2 = 40$ instances in this dataset. \\

\noindent  \textbf{ORLIB}. Large-scale problems with 1000 customers from ORLib’s facility location benchmark set, i.e., $capa$, $capb$, $capc$. They can be found at  \url{http://people.brunel.ac.uk/~mastjjb/jeb/orlib/uncapinfo.html}. The service cost $c_{ij}$ and the demand $d_i$ are given. We compute $l_{ij} = c_{ij}/d_i$. Similar to RND instances, we generate $g_{ij}$ from $[(1-\delta)\cdot l_{ij}, (1+\delta)\cdot l_{ij}]$, where $0\leq\delta< 1$.

\eat{
\begin{itemize}
\item PMPUP, a standard testbed for \textit{P-Median Problem with Users Preferences} from the well-known facility location benchmark library \textit{Discrete Location Problems} (see \url{http://www.math.nsc.ru/AP/benchmarks/Bilevel/bilevel-eng.html}).  In total, there are 30 instances (i.e., \textit{inst-333, inst-433, inst-533, ..., inst-3233}). The model presented in the library is the same as ours, making PMPUP the most suitable dataset for benchmarking the four solution approaches mentioned above.
\item RND, medium-scale instances that are randomly generated. Candidate facilities and customers are generated from a two-dimensional uniform distribution $[0,100]^2$. The distance between facility $j$ and customer $i$ (i.e., $l_{ij}$) is represented by the Euclidean distance. Each customer has a demand $d_i$, following a uniform distribution $[0,10]$. The cost of serving customer $i$ from facility $j$ is computed as $c_{ij} = d_i\cdot l_{ij}$. The disutility $g_{ij}$ is randomly generated from  $[(1-\delta)\cdot l_{ij}, (1+\delta)\cdot l_{ij}]$, where $0\leq\delta< 1$. That is, we allow $g_{ij}$ to deviate $100\delta\%$ from $l_{ij}$, and thus, the nearly facility is not necessarily the most attractive facility to customers. Furthermore, we consider the following combinations of problem sizes and parameters: $(|I|,|J|) = \{(200,150),(300,100),(400,150),(500,100)\}$, $P = \{5,10,15,20,30\}$ , $\delta = \{0.3,0.5\}$. Therefore, the number of instances are $5\times 4 \times 2 = 40$.
\item ORLIB, large-scale problems with 1000 customers from ORLib’s facility location benchmark set, i.e., $capa$, $capb$, $capc$. They can be found at  \url{http://people.brunel.ac.uk/~mastjjb/jeb/orlib/uncapinfo.html}. The service cost $c_{ij}$ and the demand $d_i$ are given. We compute $l_{ij} = c_{ij}/d_i$. Similar to RND instances, we generate $g_{ij}$ from $[(1-\delta)\cdot l_{ij}, (1+\delta)\cdot l_{ij}]$, where $0\leq\delta< 1$.
\end{itemize}
}


\subsubsection{Results analysis on PMPUP instances}


Table~\ref{Tab:PMPUP} reports the computational time in seconds (CPU[s]) of the four solution approaches over the 30 instances in the PMPUP dataset. For each instance, CPU for the approach performing the best (i.e., the smallest value among {\fontfamily{qcr}\selectfont{PDRM}}, {\fontfamily{qcr}\selectfont{SRM}}, {\fontfamily{qcr}\selectfont{Benders}}, and {\fontfamily{qcr}\selectfont{Benders-AS}}) is highlighted in boldface. The last two rows show the average CPU (AVG) and the average relatively improvement (ARI) of  {\fontfamily{qcr}\selectfont{Benders-AS}}  over the others, computed as
\begin{align}
\text{ARI} =  \frac{\text{AVG(Approach)} -  \text{AVG(Benders-AS)}}{ \text{AVG(Benders-AS)}} \times 100\%
\end{align}
which indicates how much {\fontfamily{qcr}\selectfont{Benders-AS}} is faster than the benchmark approach in terms of the AVG.

\begin{table}[h]
\centering
\footnotesize
\caption{Computational time in seconds (CPU[s]) for the PMPUP dataset.}
\label{Tab:PMPUP}
\begin{tabular}{lllll}
\toprule
Inst. & \multicolumn{1}{l}{\fontfamily{qcr}\selectfont{PDRM}} &  \multicolumn{1}{l}{\fontfamily{qcr}\selectfont{SRM}} &  \multicolumn{1}{l}{\fontfamily{qcr}\selectfont{Benders}} &  \multicolumn{1}{l}{\fontfamily{qcr}\selectfont{Benders-AS}}  \\
\hline
333   & 783.9                   & 451.4                    & 266.8                       & \textbf{220.0}                  \\
433   & 211.4                   & \textbf{124.5}           & 260.4                       & 198.4                           \\
533   & 436.5                   & 97.4                     & 21.4                        & \textbf{12.9}                   \\
633   & 1546.8                  & 322.8                    & 282.2                       & \textbf{195.5}                  \\
733   & 688.2                   & 219.6                    & 278.6                       & \textbf{196.7}                  \\
833   & 905.4                   & 269.6                    & 342.3                       & \textbf{211.7}                  \\
933   & 479.3                   & 207.6                    & 241.1                       & \textbf{198.1}                  \\
1033  & 918.7                   & 245.7                    & 226.6                       & \textbf{138.7}                  \\
1133  & 682.4                   & 218.9                    & 270.5                       & \textbf{170.3}                  \\
1233  & 1084.1                  & 398.4                    & 255.5                       & \textbf{201.2}                  \\
1333  & 1039.5                  & 581.6                    & 274.3                       & \textbf{196.2}                  \\
1433  & 1527.7                  & 650.0                    & 325.1                       & \textbf{207.7}                  \\
1533  & 257.5                   & 281.3                    & 196.4                       & \textbf{119.6}                  \\
1633  & 525.8                   & 365.3                    & 255.4                       & \textbf{196.9}                  \\
1733  & 734.1                   & 265.6                    & 220.9                       & \textbf{122.1}                  \\
1833  & 883.7                   & 264.8                    & 223.4                       & \textbf{131.9}                  \\
1933  & 439.9                   & 284.3                    & 337.3                       & \textbf{183.7}                  \\
2033  & 699.9                   & 280.5                    & 280.2                       & \textbf{165.1}                  \\
2133  & 933.6                   & 336.8                    & 311.5                       & \textbf{146.4}                  \\
2233  & 692.4                   & \textbf{128.7}           & 304.6                       & 144.7                           \\
2333  & 405.0                   & \textbf{127.8}           & 291.9                       & 142.7                           \\
2433  & 500.8                   & 230.3                    & 290.4                       & \textbf{154.9}                  \\
2533  & 239.6                   & \textbf{100.8}           & 195.3                       & 123.5                           \\
2633  & 291.5                   & \textbf{116.3}           & 226.6                       & 153.1                           \\
2733  & 633.4                   & 296.6                    & 331.6                       & \textbf{184.2}                  \\
2833  & 1506.0                  & 544.5                    & 292.2                       & \textbf{205.2}                  \\
2933  & 637.3                   & \textbf{112.1}           & 203.0                       & 125.2                           \\
3033  & 754.7                   & 200.4                    & 262.1                       & \textbf{143.7}                  \\
3133  & 616.5                   & 235.0                    & 306.9                       & \textbf{143.6}                  \\
3233  & 779.9                   & 257.0                    & 328.8                       & \textbf{178.9}                  \\
AVG   & 727.9                   & 273.9                    & 263.4                       & \textbf{163.8}                  \\
ARI   & 344.46\%              & 67.23\%                & 60.87\%                   & n.a.                  \\
\bottomrule
\end{tabular}
\end{table}

According to the table, {\fontfamily{qcr}\selectfont{PDRM}} is significantly slower than the others. As shown in (\ref{model:PDRM}), the  {\fontfamily{qcr}\selectfont{PDRM}} model introduces a large number of additional (dual) variables and constraints; therefore, the size of the formulation increases dramatically, which adds to the difficulty of handling the problem. By contrast, the {\fontfamily{qcr}\selectfont{SRM}} model in (\ref{model:SRM}) works in the original decision space and thus has substantially fewer variables and constraints. This explains why {\fontfamily{qcr}\selectfont{SRM}} is on average more than 2 times faster than {\fontfamily{qcr}\selectfont{PDRM}} in this dataset.

Surprisingly, the performance of {\fontfamily{qcr}\selectfont{SRM}} is comparably good when benchmarking with {\fontfamily{qcr}\selectfont{Benders}}: The AVG of  {\fontfamily{qcr}\selectfont{SRM}} is only slightly higher than that of  {\fontfamily{qcr}\selectfont{Benders}} (273.9 seconds versus 263.4 seconds). Moreover,  for these 30 instances,  {\fontfamily{qcr}\selectfont{SRM}} has lower CPU in 16 instances. Therefore, employing the standard branch-and-cut Benders decomposition cannot effectively speed up the computation.

In effect, the efficiency of the Benders-based approaches, to a large extent, depends on the speed of generating Benders cuts (i.e., the separation of Benders cuts) when a master solution $\bar{x}$ is found. In {\fontfamily{qcr}\selectfont{Benders}}, such a process involves solving the dual subproblem [DSP$_i(\bar{x})$] using external solvers, which may not be efficient since the solver requires compiling time and [DSP$_i(\bar{x})$]  itself is a large-scale LP that is not trivial to be solved by the LP algorithms. Considering that the separation of Benders cuts is typically performed a large number of times before the optimality is verified, it is within our expectation that {\fontfamily{qcr}\selectfont{Benders}} could be significantly slow down owing to the extensive efforts made on the separation.

Fortunately, Section~\ref{SS:accelerating_Benders} has provided an analytical and efficient method for the separation. In general, when we apply the analytical separation, we observe substantial performance improvement, i.e., the AVG of  {\fontfamily{qcr}\selectfont{Benders-AS}} is 67.23\% and 60.87\% shorter than  {\fontfamily{qcr}\selectfont{SRM}} and {\fontfamily{qcr}\selectfont{Benders}}, respectively (see the last row of Table~\ref{Tab:PMPUP}). Therefore, the branch-and-cut Benders decomposition should be supported by the well-designed analytical separation to be more practically powerful.

\eat{ xxx
Inst-633: obj 165. Original value in the library is 172.
Inst-1333: obj 168. Original value in the library is 169.}

\subsubsection{Results analysis on RND instances}

\begin{figure}[h]
\begin{center}
	\subfigure[Percentage of instances solved versus time.]{
      \psfig{figure=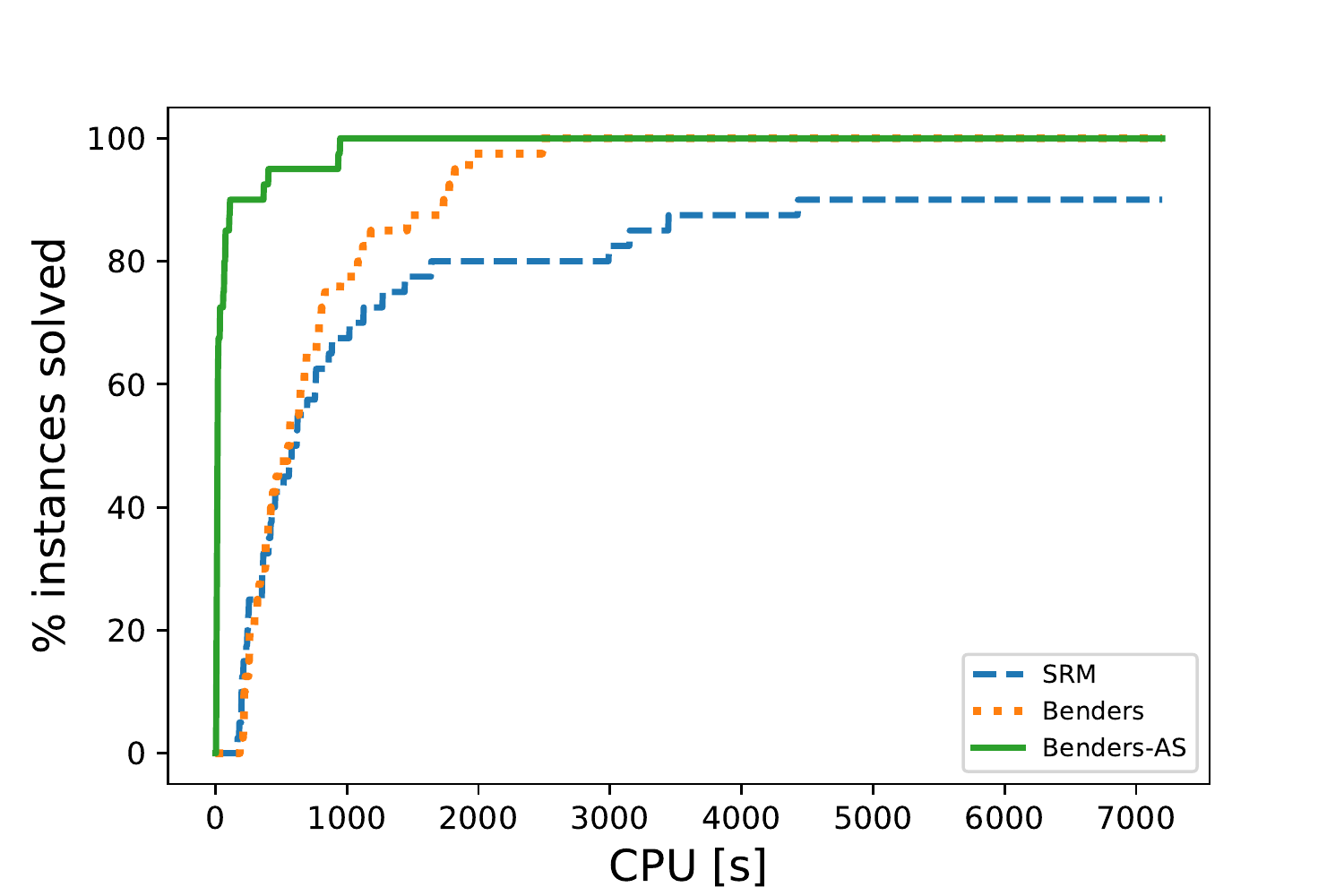,width=79mm}\label{fig:time}
    }
    \subfigure[ Box-plot  of computational time in log-scale.]{
      \psfig{figure=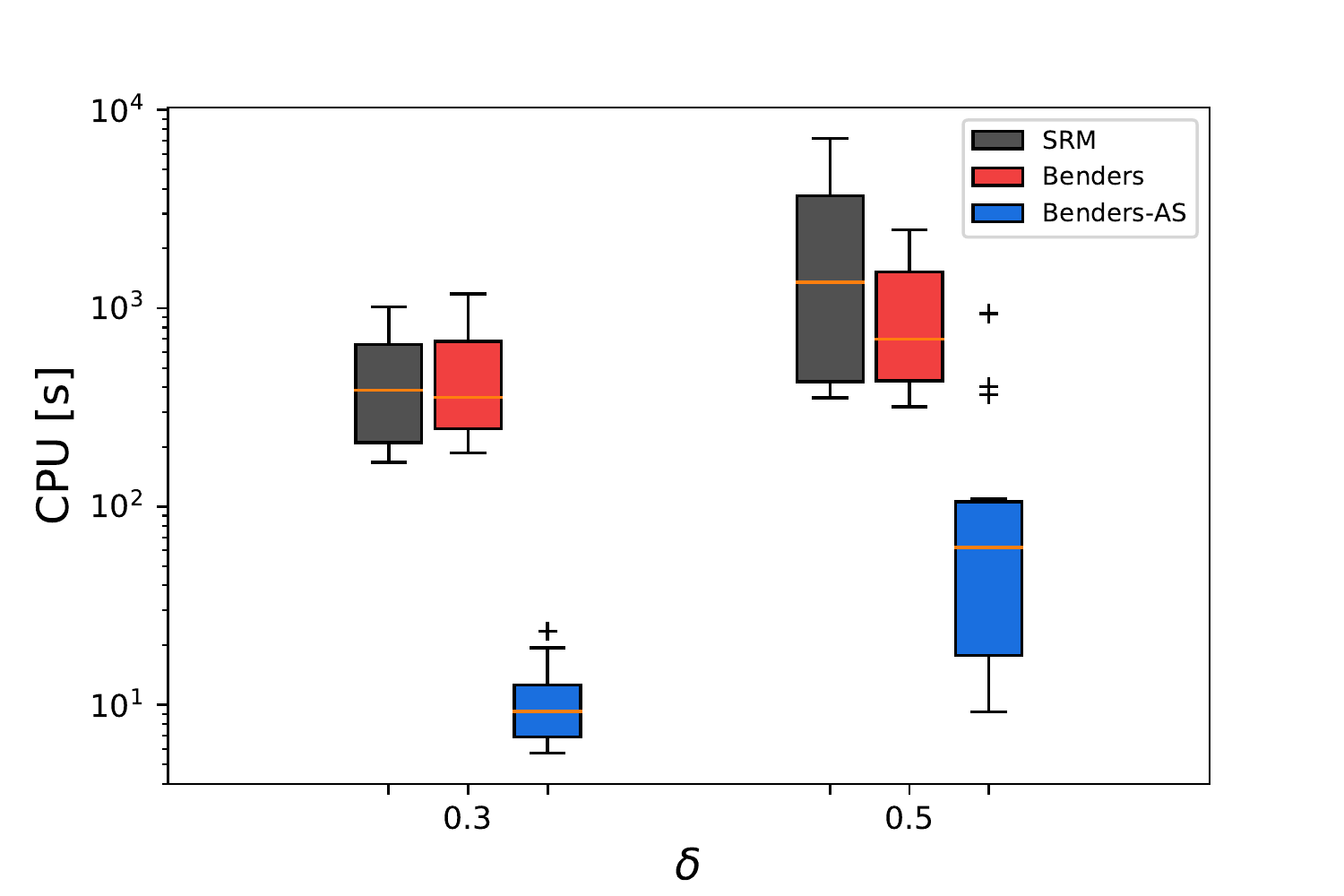,width=79mm}\label{fig:box}
      }
\end{center}
     \caption{Computational results of RND instances.}\label{fig:RND_results}
\end{figure}

We proceed to the result analysis on the RND instances. We drop {\fontfamily{qcr}\selectfont{PDRM}} since its performance is rather limited. The full computational results are provided in \ref{app:RND}. We summarize the main results in Figure~\ref{fig:RND_results}. 

Figure~\ref{fig:time} shows the percentage of instances that are solved optimally (\# instances solved) within any given CPU[s]. A point in the figure with coordinates ($m, n$) indicates that for $n$ \% of the instances, the required CPU is less than $m$ seconds to solve them to optimality. According to Figure~\ref{fig:time}, {\fontfamily{qcr}\selectfont{Benders-AS}}  outperforms the others since the green solid-line is consistently above the others, meaning that given the same CPU on the $x$-axis, {\fontfamily{qcr}\selectfont{Benders-AS}} can successfully solves more instances.  Moreover, using {\fontfamily{qcr}\selectfont{Benders-AS}}, all instances can be solved within 1000 seconds. By contrast, {\fontfamily{qcr}\selectfont{Benders}} requires roughly 2500 seconds to clear all instances, and {\fontfamily{qcr}\selectfont{SRM}} even fails in several instances within the 7200-second time limit. Therefore, in terms of the number of instances solved optimally, both Benders-based approaches outperform {\selectfont{SRM}} in this dataset.

Figure~\ref{fig:box} reports the boxplot of the computational time (in log-scale) for the above three approaches under two values of $\delta$. The line in the box indicates the median of CPU under a specific value of $\delta$. Note that a larger $\delta$ value generally indicates that customer preferences differ more from the operator's unit service costs. The difficulty of the instances certainly depends on $\delta$, and intuitively, when $\delta=0$, the model reduces to the standard facility location problem without user preferences, which can be much more efficiently solved. Thus, we expect that under a larger $\delta$ value, the instances will be more challenging. Apparently, the results in Figure~\ref{fig:RND_results}  agrees with our intuition as we observe that the instances requires longer CPU when $\delta=0.5$. Moreover, in terms of CPU, {\fontfamily{qcr}\selectfont{Benders-AS}} outperforms the others by roughly one order of magnitude, whereas the advantage of {\fontfamily{qcr}\selectfont{Benders}} over the   {\fontfamily{qcr}\selectfont{SRM}} is still not significant. This observation once again highlights the effectiveness of the proposed analytical separation in expediting the branch-and-cut Benders decomposition.

\subsubsection{Results analysis on ORLIB instances}

Our next experiment is to  investigate the performance of the approaches on the large-scale ORLIB structured data (\textit{capa, capb, capc}). Table~\ref{Tab:ORLIB} reports the results. Here, the column $rgap[\%]$ is the relative exit gap in percentages when the solution process terminates. It is computed as  $|zbb-zopt|/|zbb| \times 100$, where $zopt$ is the value of the current optimal feasible solution and $zbb$ is the best bound at termination. If an instance is solved to optimality, then $zbb = zopt$ or $rgap < 0.01\%$. The column $Ratio$ is computed by ``the CPU of an approach divided by the CPU of {\fontfamily{qcr}\selectfont{Benders-AS}}'' when the corresponding instance is solved optimally. For example, for $capa$ under  $\delta=0.1$ and $P=5$, the two $Ratio$ values under {\fontfamily{qcr}\selectfont{SRM}} and {\fontfamily{qcr}\selectfont{Benders}} are $89.8$ and $62.7$, stating that the CPUs by these two approaches are $89.8$ and $62.7$ times of the CPU by {\fontfamily{qcr}\selectfont{Benders-AS}}.  Moreover, for each problem instance, the CPU for the approach performing the best is highlighted in boldface. For those instances where all approaches cannot terminate optimally within the time limit, the smallest $rgap$ is highlighted in boldface.

According to Table~\ref{Tab:ORLIB}, we have the following observations: (i)  Consistent with the previous finding, the CPU increases with $\delta$, i.e., when the user preferences differ more from the costs of the operator, the instances becomes harder. In particular, when $\delta=0.5$ and $P=10$, the $capa$ and $capc$ instances are so changeling that none of the approaches can successfully solve them to optimality within the time limit; (ii) For all instances, {\fontfamily{qcr}\selectfont{Benders-AS}} is the best approach in this dataset. When the instance can be solved optimally, the CPU by  {\fontfamily{qcr}\selectfont{Benders-AS}} is generally more than one order of magnitude shorter. This can be directly observed from results in the column $Ratio$. Meanwhile, when the instance cannot be solved optimally within the time limit, the $rgap$ of {\fontfamily{qcr}\selectfont{Benders-AS}} is the smallest, and the value is less than 1\%, which has been small enough to guarantee the solution quality.

Based on the results of the above three experiments, we conclude that {\fontfamily{qcr}\selectfont{Benders-AS}} significantly improves the performance of the standard branch-and-cut Benders algorithm and outperforms the benchmark approaches by a large margin.

\begin{table}
\footnotesize
\centering
\caption{Computational results of the ORLIB instances.}
\label{Tab:ORLIB}
\begin{tabular}{lllllllllllll} 
\toprule
\multirow{2}{*}{Inst.} & \multicolumn{1}{c}{\multirow{2}{*}{$\delta$}} & \multicolumn{1}{c}{\multirow{2}{*}{P}} & \multicolumn{3}{c}{\fontfamily{qcr}\selectfont{SRM}}                                                          &  & \multicolumn{3}{c}{\fontfamily{qcr}\selectfont{Benders}}                                                      &  & \multicolumn{2}{c}{\fontfamily{qcr}\selectfont{Benders-AS}}                        \\ 
\cline{4-6}\cline{8-10}\cline{12-13}
                       & \multicolumn{1}{l}{}                          & \multicolumn{1}{l}{}                   & \multicolumn{1}{l}{CPU[s]} & \multicolumn{1}{c}{rgap[\%]} & \multicolumn{1}{l}{Ratio} &  & \multicolumn{1}{l}{CPU[s]} & \multicolumn{1}{l}{rgap[\%]} & \multicolumn{1}{l}{Ratio} &  & \multicolumn{1}{l}{CPU[s]} & \multicolumn{1}{l}{rgap[\%]}  \\ 
\hline
capa                   & 0.1                                           & 5                                      & 1858.3                   & 0                            & 89.8                      &  & 1298.7                   & 0                            & 62.7                      &  & \textbf{20.7 }           & 0                             \\
                       & 0.1                                           & 10                                     & 2537.0                   & 0                            & 95.0                      &  & 1813.8                   & 0                            & 67.9                      &  & \textbf{26.7 }           & 0                             \\
                       & 0.3                                           & 5                                      & 2400.9                   & 0                            & 66.3                      &  & 2226.5                   & 0                            & 61.5                      &  & \textbf{36.2 }           & 0                             \\
                       & 0.3                                           & 10                                     & 3088.6                   & 0                            & 42.6                      &  & 2516.4                   & 0                            & 34.7                      &  & \textbf{72.5 }           & 0                             \\
                       & 0.5                                           & 5                                      & 7200.0                   & 2.64                         & n.a.                      &  & 4819.1                   & 0                            & 10.7                      &  & \textbf{449.6 }          & 0                             \\
                       & 0.5                                           & 10                                     & 7200.0                   & 4.27                         & n.a.                      &  & 7200.0                   & 2.24                         & n.a.                      &  & 7200.0                   & \textbf{0.63 }                \\
                       &                                               &                                        &                          &                              &                           &  &                          &                              &                           &  &                          &                               \\
capb                   & 0.1                                           & 5                                      & 2120.1                   & 0                            & 93.0                      &  & 986.4                    & 0                            & 43.3                      &  & \textbf{22.8 }           & 0                             \\
                       & 0.1                                           & 10                                     & 1873.0                   & 0                            & 84.0                      &  & 1458.0                   & 0                            & 65.4                      &  & \textbf{22.3 }           & 0                             \\
                       & 0.3                                           & 5                                      & 2740.5                   & 0                            & 57.5                      &  & 1387.9                   & 0                            & 29.1                      &  & \textbf{47.7 }           & 0                             \\
                       & 0.3                                           & 10                                     & 2598.2                   & 0                            & 51.9                      &  & 2556.6                   & 0                            & 51.0                      &  & \textbf{50.1 }           & 0                             \\
                       & 0.5                                           & 5                                      & 7200.0                   & 1.89                         & n.a.                      &  & 5030.1                   & 0                            & 8.1                       &  & \textbf{620.5 }          & 0                             \\
                       & 0.5                                           & 10                                     & 7200.0                   & 3.51                         & n.a.                      &  & 6816.5                   & 0                            & 2.8                       &  & \textbf{2464.6 }         & 0                             \\
                       &                                               &                                        &                          &                              &                           &  &                          &                              &                           &  &                          &                               \\
capc                   & 0.1                                           & 5                                      & 1770.8                   & 0                            & 73.5                      &  & 1499.1                   & 0                            & 62.2                      &  & \textbf{24.1 }           & 0                             \\
                       & 0.1                                           & 10                                     & 2307.1                   & 0                            & 85.1                      &  & 1531.4                   & 0                            & 56.5                      &  & \textbf{27.1 }           & 0                             \\
                       & 0.3                                           & 5                                      & 2444.6                   & 0                            & 39.9                      &  & 2198.8                   & 0                            & 35.9                      &  & \textbf{61.2 }           & 0                             \\
                       & 0.3                                           & 10                                     & 2411.4                   & 0                            & 40.4                      &  & 1519.2                   & 0                            & 25.4                      &  & \textbf{59.7 }           & 0                             \\
                       & 0.5                                           & 5                                      & 7200.0                   & 4.60                         & n.a.                      &  & 7200.0                   & 0.47                         & n.a.                      &  & \textbf{1120.2 }         & 0                             \\
                       & 0.5                                           & 10                                     & 7200.0                   & 5.09                         & n.a.                      &  & 7200.0                   & 3.28                         & n.a.                      &  & 7200.0                   & \textbf{0.89 }                \\
\bottomrule
\end{tabular}
\end{table}

\subsection{Managerial implication}

In our final experiment, we investigate the importance of integrating user preferences into the facility location problem when the preference does exist. Specifically, we compare the total service costs of the operator (i.e., the objective function $\phi$) under two scenarios: (i) the operator considers user preferences and anticipates user's choices when locating facilities. This scenario is exactly the [SRM] model presented in Section~\ref{S:pd}; (ii) the operator  ignores user preferences and locates facilities based on the cost matrix $c$. In this scenario, the operator's decision is made based on the classical P-Medium problem~(see Chapter 6 of \cite{daskin2011network}).

To facilitate our discussion, we introduce the following notations. Let $x_{wt}$ be the location decision when the operator locates facilities \textit{without} considering user preferences. $x_{wt}$  can be obtained by setting $g = c$ or simply deleting Constraint~(\ref{SRM:cons4}) in the model. The operator's  \textit{actual total service cost} $\phi_{wt}$ is evaluated by holding the location decision at $x_{wt}$. Therefore,  $\phi_{wt}$ stands for the total service cost when the decision is made based on the assumption that user preferences do not exist (or are mistakenly assumed to coincide with the operator's service cost). Moreover, the \textit{optimal total service cost} $\phi$ is the optimal objective function value when the operator indeed accounts for user preferences (i.e., the best objective of [SRM]). Then, we define the relative cost increase $\Delta[\%]$ as
\begin{align}
\Delta = \frac{\phi_{wt} -\phi}{\phi} \times 100\%
\end{align}
which specifies the \textit{relative additional cost incurred when user preferences exist but is ignored by the operator when locating facilities}.

Our experiment is conducted on the PMPUP dataset since the original problem presented in the benchmark library has a similar structure to ours. Table~\ref{Tab:why} reports the operator's total service costs under three values of $P$. The first $P$ value is set to 14 since the original dataset imposes $P=14$ for all instances (except for Inst-533, marked with $*$ in the table, where $P=13$). We observe that for the original instances (i.e., the column ``$P=14$"), the average $\Delta$ is 2.44\%; therefore, ignoring user preferences will lead to an additional 2.44\% of the average service cost. Furthermore, we observe that $\Delta$ increases dramatically when $P$ increases. In particular, the average $\Delta$ under $P=30$ blows up to $84.24\%$.  This result is astonishing as it indicates that the operator will almost double its service cost; therefore, the preferences must be considered if they exist. 

\begin{table}[h]
\footnotesize
\centering
\caption{Service costs with and without the consideration of user preferences ($\phi$ and $\phi_{wt}$).}
\label{Tab:why}
\setlength{\tabcolsep}{1.5mm}{
\begin{tabular}{llllllllllll}
\toprule
\multirow{2}{*}{Inst.} & \multicolumn{3}{c}{$P=14$}                                                                  &  & \multicolumn{3}{c}{$P=20$}                                                                  &  & \multicolumn{3}{c}{$P=30$}                                                                   \\
\cline{2-4}\cline{6-8}\cline{10-12}
                       & \multicolumn{1}{l}{$\phi_{wt}$} & \multicolumn{1}{l}{$\phi$} & \multicolumn{1}{l}{$\Delta$[\%]} &~~~  &  \multicolumn{1}{l}{$\phi_{wt}$} & \multicolumn{1}{l}{$\phi$} & \multicolumn{1}{l}{$\Delta$[\%]} &~~~  & \multicolumn{1}{l}{$\phi_{wt}$} & \multicolumn{1}{l}{$\phi$} & \multicolumn{1}{l}{$\Delta$[\%]}  \\
\hline
333                    & 187                        & 172                         & 8.72                             &  & 139                        & 104                         & 33.65                            &  & 142                        & 87                          & 63.22                             \\
433                    & 156                        & 156                         & 0                                &  & 129                        & 100                         & 29.00                            &  & 134                        & 71                          & 88.73                             \\
533*                   & 188                        & 188                         & 0                                &  & 131                        & 98                          & 33.67                            &  & 133                        & 73                          & 82.19                             \\
633                    & 172                        & 165                         & 4.24                             &  & 135                        & 102                         & 32.35                            &  & 130                        & 77                          & 68.83                             \\
733                    & 160                        & 159                         & 0.63                             &  & 114                        & 91                          & 25.27                            &  & 135                        & 64                          & 110.94                            \\
833                    & 181                        & 170                         & 6.47                             &  & 142                        & 108                         & 31.48                            &  & 118                        & 81                          & 45.68                             \\
933                    & 160                        & 160                         & 0                                &  & 146                        & 102                         & 43.14                            &  & 142                        & 77                          & 84.42                             \\
1033                   & 165                        & 159                         & 3.77                             &  & 121                        & 93                          & 30.11                            &  & 128                        & 69                          & 85.51                             \\
1133                   & 174                        & 163                         & 6.75                             &  & 115                        & 98                          & 17.35                            &  & 123                        & 85                          & 44.71                             \\
1233                   & 167                        & 163                         & 2.45                             &  & 111                        & 89                          & 24.72                            &  & 109                        & 63                          & 73.02                             \\
1333                   & 169                        & 168                         & 0.60                             &  & 141                        & 102                         & 38.24                            &  & 152                        & 79                          & 92.41                             \\
1433                   & 176                        & 172                         & 2.33                             &  & 138                        & 102                         & 35.29                            &  & 132                        & 81                          & 62.96                             \\
1533                   & 152                        & 152                         & 0                                &  & 115                        & 90                          & 27.78                            &  & 108                        & 64                          & 68.75                             \\
1633                   & 156                        & 156                         & 0                                &  & 127                        & 93                          & 36.56                            &  & 134                        & 72                          & 86.11                             \\
1733                   & 160                        & 152                         & 5.26                             &  & 138                        & 94                          & 46.81                            &  & 128                        & 73                          & 75.34                             \\
1833                   & 154                        & 154                         & 0                                &  & 139                        & 91                          & 52.75                            &  & 157                        & 63                          & 149.21                            \\
1933                   & 160                        & 158                         & 1.27                             &  & 127                        & 92                          & 38.04                            &  & 127                        & 69                          & 84.06                             \\
2033                   & 161                        & 161                         & 0                                &  & 109                        & 96                          & 13.54                            &  & 144                        & 75                          & 92.00                             \\
2133                   & 176                        & 166                         & 6.02                             &  & 138                        & 108                         & 27.78                            &  & 141                        & 83                          & 69.88                             \\
2233                   & 155                        & 154                         & 0.65                             &  & 110                        & 91                          & 20.88                            &  & 155                        & 75                          & 106.67                            \\
2333                   & 160                        & 155                         & 3.23                             &  & 144                        & 98                          & 46.94                            &  & 153                        & 76                          & 101.32                            \\
2433                   & 159                        & 155                         & 2.58                             &  & 142                        & 103                         & 37.86                            &  & 155                        & 81                          & 91.36                             \\
2533                   & 147                        & 147                         & 0                                &  & 136                        & 91                          & 49.45                            &  & 119                        & 68                          & 75.00                             \\
2633                   & 156                        & 156                         & 0                                &  & 113                        & 94                          & 20.21                            &  & 146                        & 73                          & 100.00                            \\
2733                   & 173                        & 159                         & 8.81                             &  & 119                        & 102                         & 16.67                            &  & 133                        & 76                          & 75.00                             \\
2833                   & 164                        & 161                         & 1.86                             &  & 132                        & 95                          & 38.95                            &  & 115                        & 73                          & 57.53                             \\
2933                   & 152                        & 152                         & 0                                &  & 117                        & 91                          & 28.57                            &  & 119                        & 65                          & 83.08                             \\
3033                   & 162                        & 157                         & 3.18                             &  & 120                        & 93                          & 29.03                            &  & 105                        & 62                          & 69.35                             \\
3133                   & 155                        & 155                         & 0                                &  & 109                        & 87                          & 25.29                            &  & 165                        & 63                          & 161.90                            \\
3233                   & 162                        & 155                         & 4.52                             &  & 128                        & 97                          & 31.96                            &  & 137                        & 77                          & 77.92                             \\
\textbf{AVG}           & \textbf{164.0}            & \textbf{159.1}             & \textbf{2.24}                    &  & \textbf{127.5}            & \textbf{96.5}              & \textbf{32.1}                   &  & \textbf{134.0}            & \textbf{73.2}              & \textbf{84.24}                    \\
\bottomrule
\end{tabular}}
\end{table}

\begin{figure}[h]
\begin{center}
	\subfigure[Inst-333]{
      \psfig{figure=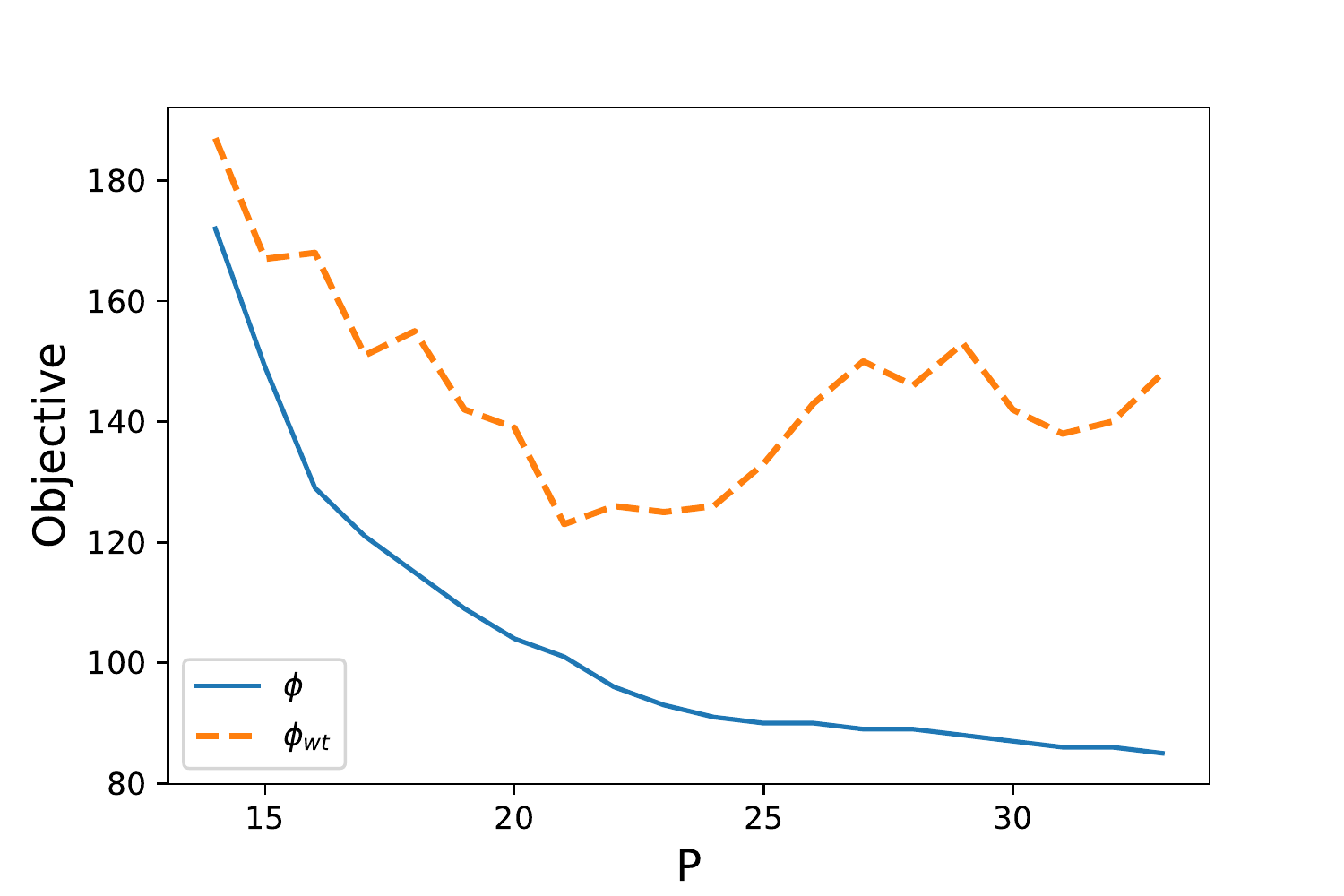,width=78mm}\label{fig:333}
    }
    \subfigure[Inst-1333]{
      \psfig{figure=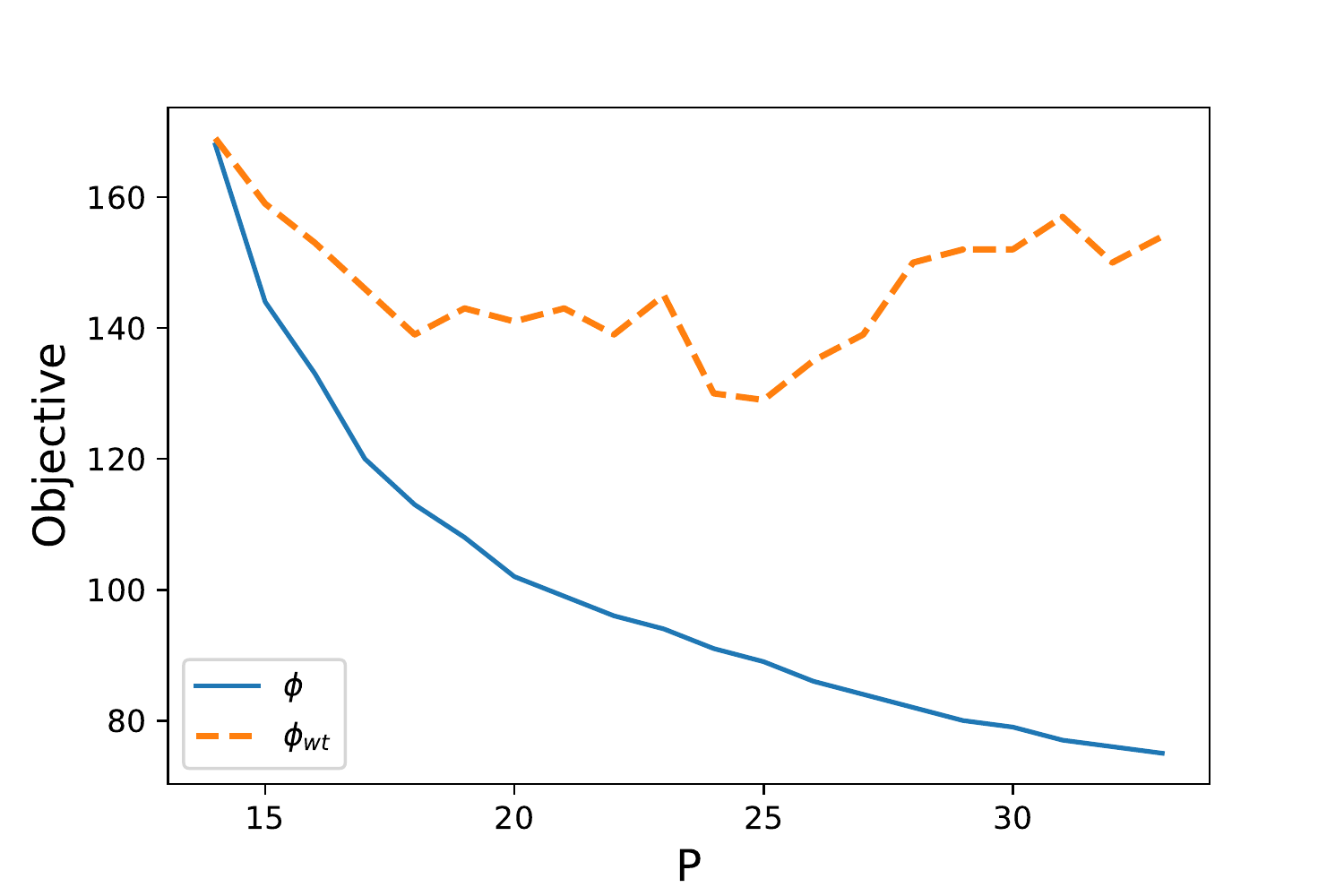,width=78mm}\label{fig:1333}
      }
\end{center}
     \caption{Service costs with and without the consideration of user preferences versus the number of open facilities.}\label{fig:implication}
\end{figure}

Another interesting observation is that when user preferences are ignored, opening more facilities may lead to a higher service cost. For example, in the Inst-333, when the operator opens 20 facilities, $\phi_{wt}$ is 139; however, when 10 more facilities are open, $\phi_{wt}$ increases to 142. To give a direct view, we plot how $\phi$ and $\phi_{wt}$ evolve with $P$ in Figure~\ref{fig:implication}. The figure shows $\phi_{wt}$ could increase when $P$ increases. On the contrary, we can clearly observe the decreasing trend of $\phi$, stating that when the preferences are considered in the operator's decision stage, opening more facilities will indeed reduce the service cost.

To summarize, the above experiment demonstrates that the operator must take user preferences into consideration and correctly anticipate user's choices; otherwise, the operator could suffer from a substantially higher cost, and opening more facilities could unexpectedly result in additional service costs as well. These results further justify the usefulness of the P-median problem with user preferences.

\section{Conclusion}
\label{S:cl}

This paper studied the exact solution approach for the P-median facility location problem with user preferences (PUP). By exploring the problem structure, we proved that in a CAC-based MILP model~\citep{casas2017solving}, the high-dimensional binary variables used to model customer preferences and facility choices can be relaxed to continuous variables. Based on this, we proposed a branch-and-cut algorithm where Benders separation (i.e., the procedure of generating Benders cuts) at integer nodes of the searching tree was performed leveraging external solvers. However, such a standard Benders approach was  not efficient enough. Therefore, we further proposed an acceleration technique to enhance the algorithm performance. Using a broad testbed, our computational experiments demonstrated that the proposed algorithm outperformed several benchmark approaches by a large margin and was able to handle large-scale instances satisfactorily.  We also conducted sensitivity analysis and observed that when user preferences indeed exist, the operator must consider customer preferences and correctly anticipate the choices to avoid an unnecessarily high service provision cost. Furthermore, ignoring the preferences may lead to an ironic situation where opening more facilities could result in additional service costs.

There are limitations and potential future research directions. Firstly, as shown in~\cite{fischetti2017redesigning}, one may further enhance the Benders decomposition by designing a proper cut loop strategy and stabilization at the root node.  However, fractional location solutions arise when stabilization is implemented. Noting that the analytical separation is only applicable to integer nodes, we are thus unable to generate Benders cuts efficiently when stabilizing the solution, and an alternative method must be used instead. This means that an efficient cut loop is not ready to been obtained, and its effectiveness on the algorithm also needs more careful investigations. Secondly, the [SRM] formulation leveraged a \textit{closest assignment constraint} (CAC) to transform the bilevel model into a MILP. The whole Benders approach was built upon this formulation. It is then interesting to explore other CACs since their strengths and sizes are different~\citep{espejo2012closest}. It is possible to have a powerful Benders algorithm that is developed based on other CACs. Finally, since the operator needs to forecast customer preferences to anticipate the choices, there exists the possibility that the forecasted preference is subject to estimation errors. In this case, we should build a stochastic model or a robust optimization model, and we would like to leave the model formulation and the algorithm development for future research.

\appendix

\section{Primal-dual reformulation model}
\label{app:pdr}

This appendix presents a single-level reformulation model for [PUP], which relies on the primal-dual optimality conditions of the lower-level problem. Similar approaches can be found in \cite{casas2017solving},  \cite{casas2018approximating}. Here, we briefly describe the reformulation procedure. Given the operator's location decision $x$,  the lower-level problem is
\begin{subequations}\label{pdr:sub:}
\begin{alignat}{10}
\min_{y \in \R_+} &\sum_{i \in I}\sum_{j \in J} \pi_{ij}y_{ij} \\
\text{st.}~&\sum_{j \in J} y_{ij} = 1  \quad \forall i \in I         & (\alpha_i)\\
&y_{ij} \leq x_j  \quad \forall i \in I, j \in J   \qquad & (\beta_{ij})
\end{alignat}
\end{subequations}
The $\alpha_i$ and $\beta_{ij}$ in the parentheses are the dual variables associated with the constraints. We can then rewrite the lower-level problem with its KKT conditions, i.e., 
\begin{subequations}\label{pdr:KKT}
\begin{alignat}{10}
&\sum_{j \in J} y_{ij} = 1   & \forall i \in I \\
& y_{ij} \leq x_j & \forall i \in I, j \in J \\
&y_{ij} \geq 0  & \forall i \in I, j \in J \\
& \alpha_i + \beta_{ij} \leq g_{ij}  & \forall i \in I, j \in J \\
\label{prd:bilinear1}& y_{ij} ( \alpha_i + \beta_{ij} - g_{ij}) = 0  \qquad & \forall i \in I, j \in J \\
\label{prd:bilinear2}& \beta_{ij} (y_{ij} - x_j ) = 0  & \forall i \in I, j \in J \\
& \beta_{ij} \leq 0  & \forall i \in I, j \in J
\end{alignat}
\end{subequations}
where (\ref{prd:bilinear1}) and (\ref{prd:bilinear2}) are bilinear functions. Noting both $y_{ij}$ and $x_j$ will take 0/1 in the optimal solution and the maximum value of the matrix $\pi$ is 1 by definition, these bilinear functions can exactly linearized by 
\begin{subequations}\label{pdr:linearized}
\begin{alignat}{2}
\label{prd:linear1}&\alpha_i + \beta_{ij} - \pi_{ij} \geq - (1-y_{ij}) \qquad & \forall i \in I, j \in J \\
\label{prd:linear2}& \beta_{ij} \geq - (1 +  y_{ij} - x_j )& \forall i \in I, j \in J 
\end{alignat}
\end{subequations}
Then, [PUP] is equivalent to
\begin{subequations}\label{model:PDRM}
\begin{alignat}{10}
 \min~& \sum_{i \in I} \sum_{j \in J} c_{ij}y_{ij} \\
\text{st.}~ & (x,y,\alpha,\beta) \in (\ref{pdr:KKT}) \\
\textbf{[PDRM]} \qquad &\sum_{j \in J} x_{j} = P \\
& x_{j} \in\{0,1\} \quad \forall j \in J
\end{alignat}
\end{subequations}
which is referred to as the primal-dual-reformulation model. In our preliminary computational experiment, we observe that setting $y_{ij} \in \{0,1\}$ (which is feasible and equivalent) makes Gurobi run faster and enhances the numerical stability as well.

\section{Full computational results of RND instances}
\label{app:RND}

This appendix presents the computational results of the RND instances under the 2-hour time limit in Table~\ref{Tab:RND}. The RND dataset is generated with fixed random seed in Python files and can be provided upon request.

\begin{table}[htbp]
\footnotesize
\centering
\caption{Computational results of the RND instances.}
\label{Tab:RND}
\setlength{\tabcolsep}{1.8mm}{
\begin{tabular}{llllllllllllll} 
\toprule
\multirow{2}{*}{$\delta$} & \multicolumn{1}{c}{\multirow{2}{*}{$(|I|,|J|)$}} & \multicolumn{1}{c}{\multirow{2}{*}{P~~}} & \multicolumn{3}{c}{\fontfamily{qcr}\selectfont{SRM}}                                                              &  & \multicolumn{3}{c}{\fontfamily{qcr}\selectfont{Benders}}                                                          & ~~ & \multicolumn{3}{c}{\fontfamily{qcr}\selectfont{Benders-AS}}                                                        \\ 
\cline{4-6}\cline{8-10}\cline{12-14}
                          & ~~                         & \multicolumn{1}{c}{}                   & \multicolumn{1}{c}{CPU[s]} & \multicolumn{1}{c}{rgap[\%]} & \multicolumn{1}{c}{$\phi$} &  & \multicolumn{1}{c}{CPU[s]} & \multicolumn{1}{c}{rgap[\%]} & \multicolumn{1}{c}{$\phi$} &  & \multicolumn{1}{c}{CPU[s]} & \multicolumn{1}{c}{rgap[\%]} & \multicolumn{1}{c}{$\phi$}  \\ 
\hline
0.3                       & (200,150)                                     & 5                                      & 201.5                    & 0                            & 19883                      &  & 226.1                    & 0                            & 19883                      &  & 6.3                      & 0                            & 19883                       \\
                          &                                               & 10                                     & 212.4                    & 0                            & 13539                      &  & 296.1                    & 0                            & 13539                      &  & 8.3                      & 0                            & 13539                       \\
                          &                                               & 15                                     & 196.2                    & 0                            & 10701                      &  & 186.9                    & 0                            & 10701                      &  & 6.9                      & 0                            & 10701                       \\
                          &                                               & 20                                     & 241.1                    & 0                            & 9215                       &  & 211.5                    & 0                            & 9215                       &  & 6.9                      & 0                            & 9215                        \\
                          &                                               & 30                                     & 180.1                    & 0                            & 7506                       &  & 217.6                    & 0                            & 7506                       &  & 5.7                      & 0                            & 7506                        \\
                          & (300,100)                                     & 5                                      & 233.9                    & 0                            & 29993                      &  & 257.1                    & 0                            & 29993                      &  & 8.4                      & 0                            & 29993                       \\
                          &                                               & 10                                     & 253.9                    & 0                            & 20404                      &  & 217.7                    & 0                            & 20404                      &  & 7.1                      & 0                            & 20404                       \\
                          &                                               & 15                                     & 166.8                    & 0                            & 16523                      &  & 253.0                    & 0                            & 16523                      &  & 7.3                      & 0                            & 16523                       \\
                          &                                               & 20                                     & 247.4                    & 0                            & 14479                      &  & 330.3                    & 0                            & 14479                      &  & 6.7                      & 0                            & 14479                       \\
                          &                                               & 30                                     & 198.1                    & 0                            & 12155                      &  & 257.5                    & 0                            & 12155                      &  & 5.7                      & 0                            & 12155                       \\
                          & (400,150)                                     & 5                                      & 757.5                    & 0                            & 39161                      &  & 676.2                    & 0                            & 39161                      &  & 18.4                     & 0                            & 39161                       \\
                          &                                               & 10                                     & 885.1                    & 0                            & 27180                      &  & 768.3                    & 0                            & 27180                      &  & 23.5                     & 0                            & 27180                       \\
                          &                                               & 15                                     & 1018.3                   & 0                            & 21666                      &  & 568.9                    & 0                            & 21666                      &  & 19.4                     & 0                            & 21666                       \\
                          &                                               & 20                                     & 763.9                    & 0                            & 18494                      &  & 1179.5                   & 0                            & 18494                      &  & 15.6                     & 0                            & 18494                       \\
                          &                                               & 30                                     & 860.8                    & 0                            & 15138                      &  & 786.3                    & 0                            & 15138                      &  & 12.5                     & 0                            & 15138                       \\
                          & (500,100)                                     & 5                                      & 623.5                    & 0                            & 49530                      &  & 635.6                    & 0                            & 49530                      &  & 12.1                     & 0                            & 49530                       \\
                          &                                               & 10                                     & 618.2                    & 0                            & 33744                      &  & 382.1                    & 0                            & 33744                      &  & 11.6                     & 0                            & 33744                       \\
                          &                                               & 15                                     & 578.1                    & 0                            & 27148                      &  & 691.3                    & 0                            & 27148                      &  & 13.0                     & 0                            & 27148                       \\
                          &                                               & 20                                     & 517.2                    & 0                            & 23445                      &  & 646.4                    & 0                            & 23445                      &  & 10.8                     & 0                            & 23445                       \\
                          &                                               & 30                                     & 559.5                    & 0                            & 19686                      &  & 805.8                    & 0                            & 19686                      &  & 10.1                     & 0                            & 19686                       \\
                          &                                               &                                        &                          &                              &                            &  &                          &                              &                            &  &                          &                              &                             \\
0.5                       & (200,150)                                     & 5                                      & 418.1                    & 0                            & 20494                      &  & 566.7                    & 0                            & 20494                      &  & 15.9                     & 0                            & 20494                       \\
                          &                                               & 10                                     & 1271.0                   & 0                            & 14120                      &  & 549.6                    & 0                            & 14120                      &  & 74.8                     & 0                            & 14120                       \\
                          &                                               & 15                                     & 1125.9                   & 0                            & 11208                      &  & 419.4                    & 0                            & 11208                      &  & 33.3                     & 0                            & 11208                       \\
                          &                                               & 20                                     & 404.3                    & 0                            & 9536                       &  & 380.4                    & 0                            & 9536                       &  & 14.1                     & 0                            & 9536                        \\
                          &                                               & 30                                     & 357.0                    & 0                            & 7775                       &  & 459.4                    & 0                            & 7775                       &  & 9.2                      & 0                            & 7775                        \\
                          & (300,100)                                     & 5                                      & 362.7                    & 0                            & 31004                      &  & 434.2                    & 0                            & 31004                      &  & 21.0                     & 0                            & 31004                       \\
                          &                                               & 10                                     & 453.7                    & 0                            & 21133                      &  & 317.8                    & 0                            & 21133                      &  & 18.4                     & 0                            & 21133                       \\
                          &                                               & 15                                     & 696.6                    & 0                            & 17225                      &  & 400.4                    & 0                            & 17225                      &  & 33.1                     & 0                            & 17225                       \\
                          &                                               & 20                                     & 428.0                    & 0                            & 15005                      &  & 479.1                    & 0                            & 15005                      &  & 15.9                     & 0                            & 15005                       \\
                          &                                               & 30                                     & 353.5                    & 0                            & 12601                      &  & 342.8                    & 0                            & 12601                      &  & 11.5                     & 0                            & 12601                       \\
                          & (400,150)                                     & 5                                      & 3149.0                   & 0                            & 40361                      &  & 1818.4                   & 0                            & 40361                      &  & 109.3                    & 0                            & 40361                       \\
                          &                                               & 10                                     & 7200.0                   & 1.48                         & 28185                      &  & 1928.2                   & 0                            & 28185                      &  & 366.9                    & 0                            & 28185                       \\
                          &                                               & 15                                     & 7200.0                   & 2.01                         & 22685                      &  & 2485.3                   & 0                            & 22685                      &  & 946.9                    & 0                            & 22685                       \\
                          &                                               & 20                                     & 7200.0                   & 0.52                         & 19234                      &  & 1733.7                   & 0                            & 19234                      &  & 401.9                    & 0                            & 19234                       \\
                          &                                               & 30                                     & 7200.0                   & 0.35                         & 15858                      &  & 1776.9                   & 0                            & 15858                      &  & 933.5                    & 0                            & 15858                       \\
                          & (500,100)                                     & 5                                      & 1439.6                   & 0                            & 51628                      &  & 1459.0                   & 0                            & 51628                      &  & 67.5                     & 0                            & 51628                       \\
                          &                                               & 10                                     & 4427.7                   & 0                            & 35360                      &  & 948.8                    & 0                            & 35360                      &  & 104.3                    & 0                            & 35360                       \\
                          &                                               & 15                                     & 1639.4                   & 0                            & 28220                      &  & 830.7                    & 0                            & 28220                      &  & 65.8                     & 0                            & 28220                       \\
                          &                                               & 20                                     & 3443.4                   & 0                            & 24631                      &  & 1081.8                   & 0                            & 24631                      &  & 58.6                     & 0                            & 24631                       \\
                          &                                               & 30                                     & 2989.9                   & 0                            & 20807                      &  & 1119.4                   & 0                            & 20807                      &  & 76.2                     & 0                            & 20807                       \\
\bottomrule
\end{tabular}}
\end{table}

\eat{

\section{Proof of Lemma~\ref{lemma_relax}}
\label{app:lemma_relax}
Suppose a location decision is made at $\bar{x}$. Let $\bar{\tau}$ be the set of open facilities, i.e., $\bar{\tau} = \left\{ \forall j \in J \mid \bar{x}_j = 1 \right\}$. Based on the condition $y_{ij} \leq \bar{x}_j$, Constraint~(\ref{SRM:cons4}) can be restated as $ \sum_{j \in \bar{\tau} }\pi_{ij}y_{ij} \leq \min_{j \in \bar{\tau}} \pi_{ij}, \forall i \in I$. For notation simplicity, we drop subscript $i$ in this proof. 

Let $m$ denote the most preferred open facility for customers, i.e., $m =  \arg \min_{j\in \bar{\tau}} \pi_{j}$. We have $ \sum_{j \in \bar{\tau} }\pi_{j}y_{j} \leq \pi_{m}$. Note that $\pi_m < \pi_j, \forall j \in  \bar{\tau} \setminus \{m\}$, and $\sum_{j \in \bar{\tau}} y_j = 1$. Therefore, $ \sum_{j \in \bar{\tau} }\pi_{j}y_{j} \leq \pi_{m}$ holds only if $y_m = 1$ and $y_j = 0, \forall j \in \bar{\tau} \setminus \{m\}$.  This indicates that relaxing $y_j \in \{0,1\}$ to $ y_j \geq 0$ will not change the solution since $y_j$ in [SRM] will automatically be an integer.
xxx
}

\eat{
\section{Reformulation based on the approach in VK2010}

\begin{align}
T_{ij} = \left\{ \forall k \in J \mid g_{ik} > g_{ij} \right\}, \forall i \in I, j \in J
\end{align}

Set of facilities that are dominated by facility $j$ for customer $i$.

\begin{align}
\sum_{k \in T_{ij}}y_{ik} \leq 1 - x_j \quad \forall i \in I, j \in J 
\end{align}

\begin{subequations}\label{model:SRM}
\begin{alignat}{10}
\label{VK:cons1} \min~& \sum_{i \in I} \sum_{j \in J} c_{ij}y_{ij} \\
\label{VK:cons2}\text{st.}~& \sum_{j \in J} y_{ij} = 1 \quad \forall i \in I \\
\label{VK:cons3}\textbf{[VK]} \qquad \quad  & y_{ij} \leq x_j  \quad \forall i \in I, j \in J \\
\label{VK:cons4}& \sum_{k \in T_{ij}}y_{ik} \leq 1 - x_j \quad \forall i \in I, j \in J \\
\label{VK:cons5}& x_{j} \in \Omega \quad \forall j \in J\\
\label{VK:cons6}& y_{ij} \in \{0,1\} \quad \forall i \in I, j \in J
\end{alignat}
\end{subequations}
}

\eat{

\section*{Acknowledgement}
The authors would like to thank the editor Jos{\'e} Fernando Oliveira and five anonymous reviewers for their detailed and constructive comments that have significantly enhanced the quality of this work.

\subsection{RND instances 1}

facilities and customers: uniform $[0,100]^2$. Euclidean distance $l$.

Each customer demand:  uniform $[0,10]^2$

Set $c = d\cdot l$

Set $g = l + \epsilon$, where $\epsilon \in [0,\Delta]$.

$P = \{5,6,7,8,9,10\}$

$\Delta = \{10,30\}$

$(I,J) = \{(200,150),(300,100),(400,150),(500,100)\}$

 Time limit is 7200 seconds.

\begin{table}[h]
\small
\centering
\caption{Computational results of the RND instances. Each row shows the average statistics over 6 values of $P$.}
\label{Tab:Result200-05}
\setlength{\tabcolsep}{1.5mm}{
\begin{tabular}{rrrrrrrrrrrrrr}
\toprule
\multirow{2}{*}{$\Delta$} & \multicolumn{1}{c}{\multirow{2}{*}{$|I|$}} & \multicolumn{1}{c}{\multirow{2}{*}{$|J|$}} & \multicolumn{3}{c}{MILP}                                                          & ~~ & \multicolumn{3}{c}{Benders}                                                       &  ~~& \multicolumn{3}{c}{Benders-AS}                                                     \\
\cline{4-6}\cline{8-10}\cline{12-14}
                           & \multicolumn{1}{c}{}                       & \multicolumn{1}{c}{}                       & \multicolumn{1}{c}{t[s]} & \multicolumn{1}{c}{rgap[\%]} & \multicolumn{1}{c}{Obj} &  & \multicolumn{1}{c}{t[s]} & \multicolumn{1}{c}{rgap[\%]} & \multicolumn{1}{c}{Obj} &  & \multicolumn{1}{c}{t[s]} & \multicolumn{1}{c}{rgap[\%]} & \multicolumn{1}{c}{Obj}  \\
\hline
10                         & 200                                        & 150                                        & 176.4                    & 0                            & 16175                 &  & 258.7                    & 0                            & 16175                 &  & 6.1                      & 0                            & 16175                  \\
                           & 300                                        & 100                                        & 191.3                    & 0                            & 24418                 &  & 194.9                    & 0                            & 24418                 &  & 6.1                      & 0                            &24418                  \\
                           & 400                                        & 150                                        & 647.3                    & 0                            & 32003                 &  & 737.1                    & 0                            & 32003                 &  & 14.3                     & 0                            & 32003                  \\
                           & 500                                        & 100                                        & 505.3                    & 0                            & 40335                 &  & 485.5                    & 0                            & 40335                 &  & 9.8                      & 0                            & 40335                  \\
                           &                                            &                                            &                          &                              &                         &  &                          &                              &                         &  &                          &                              &                          \\
30                         & 200                                        & 150                                        & 2235.8                   & 0.22                         & 17023                 &  & 877.8                    & 0                            & 17023                 &  & 115.9                    & 0                            & 17023                  \\
                           & 300                                        & 100                                        & 1917.1                   & 0                            & 25816                 &  & 583.8                    & 0                            & 25816                 &  & 80.4                     & 0                            & 25816                  \\
                           & 400                                        & 150                                        & 6321.8                   & 2.51                         & 34112                 &  & 4090.8                   & 0.94                         & 34112                 &  & 2708.0                   & 0.78                         & 34109                  \\
                           & 500                                        & 100                                        & 4812.8                   & 1.57                         & 43121                 &  & 1980.9                   & 0                            & 43121                 &  & 602.0                    & 0                            & 43121                  \\
\bottomrule
\end{tabular}}
\end{table}

\begin{figure}[h]
\begin{center}
      \psfig{figure=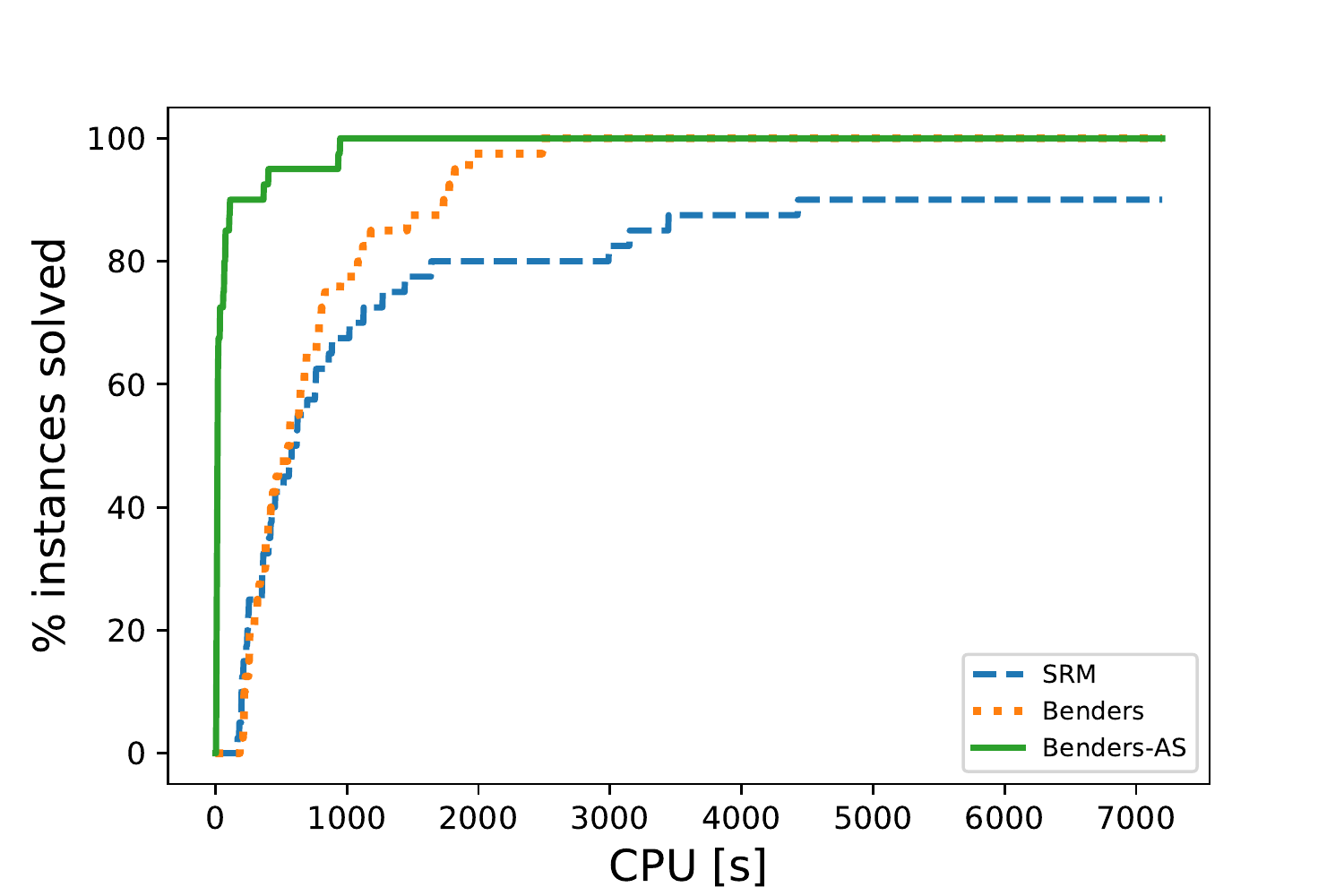,width=100mm}\label{fig:instances_solved}
\end{center}
\caption{RND instances}
\end{figure}

\eat{

\begin{table}[h]
\footnotesize
\centering
\caption{Original DLPPUP instances where $P = 13$ or 14 . Increase in the cost when the preference is not considered.}
\label{Tab:why}
\setlength{\tabcolsep}{0.5mm}{
\begin{tabular}{rrrrrrrrrrr}
\toprule
\multirow{2}{*}{Inst.} & \multicolumn{4}{c}{Objective}                                                                                               &~~~~  & \multirow{2}{*}{Inst.} & \multicolumn{4}{c}{Objective~}                                                                                               \\
\cline{2-5}\cline{8-11}
                       & \multicolumn{1}{c}{Predicted} & \multicolumn{1}{c}{Actual} & \multicolumn{1}{c}{Optimal} & \multicolumn{1}{c}{Rel-Loss[\%]} &  &                        & \multicolumn{1}{c}{Predicted} & \multicolumn{1}{c}{Actual} & \multicolumn{1}{c}{Optimal} & \multicolumn{1}{c}{Rel-Loss[\%]}  \\
\hline
333                    & 147                           & 187                        & 172                         & 8.72                             &  & 1833                   & 139                           & 154                        & 154                         & 0                                 \\
433                    & 145                           & 156                        & 156                         & 0                                &  & 1933                   & 137                           & 160                        & 158                         & 1.27                              \\
533                    & 177                           & 188                        & 188                         & 0                                &  & 2033                   & 140                           & 161                        & 161                         & 0                                 \\
633                    & 144                           & 172                        & 165                         & 4.24                             &  & 2133                   & 138                           & 176                        & 166                         & 6.02                              \\
733                    & 137                           & 160                        & 159                         & 0.63                             &  & 2233                   & 121                           & 155                        & 154                         & 0.65                              \\
833                    & 144                           & 181                        & 170                         & 6.47                             &  & 2333                   & 133                           & 160                        & 155                         & 3.23                              \\
933                    & 130                           & 160                        & 160                         & 0                                &  & 2433                   & 139                           & 159                        & 155                         & 2.58                              \\
1033                   & 138                           & 165                        & 159                         & 3.77                             &  & 2533                   & 131                           & 147                        & 147                         & 0                                 \\
1133                   & 147                           & 174                        & 163                         & 6.75                             &  & 2633                   & 132                           & 156                        & 156                         & 0                                 \\
1233                   & 142                           & 167                        & 163                         & 2.45                             &  & 2733                   & 139                           & 173                        & 159                         & 8.81                              \\
1333                   & 140                           & 169                        & 168                         & 0.60                             &  & 2833                   & 137                           & 164                        & 161                         & 1.86                              \\
1433                   & 152                           & 176                        & 172                         & 2.33                             &  & 2933                   & 124                           & 152                        & 152                         & 0                                 \\
1533                   & 133                           & 152                        & 152                         & 0                                &  & 3033                   & 137                           & 162                        & 157                         & 3.18                              \\
1633                   & 141                           & 156                        & 156                         & 0                                &  & 3133                   & 141                           & 155                        & 155                         & 0                                 \\
1733                   & 134                           & 160                        & 152                         & 5.26                             &  & 3233                   & 129                           & 162                        & 155                         & 4.52                              \\
\bottomrule
\end{tabular}}
\end{table}
}

\eat{
\begin{table}
\centering
\caption{Computational results for the PMPUP dataset.}
\label{Tab:PMPUP}
\begin{tabular}{rrrrrrrrrr}
\hline
\multirow{2}{*}{Inst} & \multicolumn{2}{c}{MILP}                           &  & \multicolumn{2}{c}{Benders}                            &  & \multicolumn{3}{c}{Bender+AS}                                                 \\
\cline{2-3}\cline{5-6}\cline{8-10}
                          & \multicolumn{1}{c}{t[s]} & \multicolumn{1}{c}{\#N} &  & t[s]                      & \#N                        &  & \multicolumn{1}{c}{t[s]} & \multicolumn{1}{c}{\#N} & \multicolumn{1}{c}{Obj}  \\
\hline
333                       & 451.4                    & 313349                  &  & 266.8                     & 760994                     &  & \textbf{220.0}           & 755341                  & 172                      \\
433                       & \textbf{124.5}           & 83255                   &  & 260.4                     & 552144                     &  & 198.4                    & 665937                  & 156                      \\
533                       & 97.4                     & 53587                   &  & 21.4                      & 61399                      &  & \textbf{12.9}            & 62950                   & 188                      \\
633                       & 322.8                    & 219116                  &  & 282.2                     & 723295                     &  & \textbf{195.5}           & 730947                  & $^*165$                     \\
733                       & 219.6                    & 184225                  &  & 278.6                     & 739050                     &  & \textbf{196.7}           & 698104                  & 159                      \\
833                       & 269.6                    & 201037                  &  & 342.3                     & 738552                     &  & \textbf{211.7}           & 719946                  & 170                      \\
933                       & 207.6                    & 141560                  &  & 241.1                     & 622913                     &  & \textbf{198.1}           & 716351                  & 160                      \\
1033                      & 245.7                    & 175062                  &  & 226.6                     & 623777                     &  & \textbf{138.7}           & 536114                  & 159                      \\
1133                      & 218.9                    & 159760                  &  & 270.5                     & 704502                     &  & \textbf{170.3}           & 649630                  & 163                      \\
1233                      & 398.4                    & 268616                  &  & 255.5                     & 580285                     &  & \textbf{201.2}           & 712997                  & 163                      \\
1333                      & 581.6                    & 393798                  &  & 274.3                     & 743414                     &  & \textbf{196.2}           & 759721                  & $^*168$                     \\
1433                      & 650.0                    & 428436                  &  & 325.1                     & 797505                     &  & \textbf{207.7}           & 809408                  & 172                      \\
1533                      & 281.3                    & 185474                  &  & 196.4                     & 563567                     &  & \textbf{119.6}           & 529596                  & 152                      \\
1633                      & 365.3                    & 232346                  &  & 255.4                     & 613526                     &  & \textbf{196.9}           & 696730                  & 156                      \\
1733                      & 265.6                    & 197428                  &  & 220.9                     & 497231                     &  & \textbf{122.1}           & 452213                  & 152                      \\
1833                      & 264.8                    & 196625                  &  & 223.4                     & 635657                     &  & \textbf{131.9}           & 565629                  & 154                      \\
1933                      & 284.3                    & 217233                  &  & 337.3                     & 635825                     &  & \textbf{183.7}           & 616006                  & 158                      \\
2033                      & 280.5                    & 220177                  &  & 280.2                     & 657664                     &  & \textbf{165.1}           & 654264                  & 161                      \\
2133                      & 336.8                    & 217067                  &  & 311.5                     & 706950                     &  & \textbf{146.4}           & 571052                  & 166                      \\
2233                      & \textbf{128.7}           & 89847                   &  & 304.6                     & 654347                     &  & 144.7                    & 550327                  & 154                      \\
2333                      & \textbf{127.8}           & 87620                   &  & 291.9                     & 669038                     &  & 142.7                    & 560988                  & 155                      \\
2433                      & 230.3                    & 135437                  &  & 290.4                     & 603530                     &  & \textbf{154.9}           & 619764                  & 155                      \\
2533                      & \textbf{100.8}           & 57386                   &  & 195.3                     & 470600                     &  & 123.5                    & 454579                  & 147                      \\
2633                      & \textbf{116.3}           & 79565                   &  & 226.6                     & 670333                     &  & 153.1                    & 594557                  & 156                      \\
2733                      & 296.6                    & 197553                  &  & 331.6                     & 727162                     &  & \textbf{184.2}           & 647388                  & 159                      \\
2833                      & 544.5                    & 389994                  &  & 292.2                     & 771520                     &  & \textbf{205.2}           & 746886                  & 161                      \\
2933                      & \textbf{112.1}           & 75945                   &  & 203.0                     & 527056                     &  & 125.2                    & 471001                  & 152                      \\
3033                      & 200.4                    & 155548                  &  & 262.1                     & 639917                     &  & \textbf{143.7}           & 630984                  & 157                      \\
3133                      & 235.0                    & 145446                  &  & 306.9                     & 692432                     &  & \textbf{143.6}           & 595924                  & 155                      \\
3233                      & 257.0                    & 215603                  &  & \multicolumn{1}{l}{328.8} & \multicolumn{1}{l}{706074} &  & \textbf{178.9}           & 640627                  & 155                      \\
Avg                   & 273.9                    &                         &  & 263.4                     &                            &  & \textbf{163.8}           &                         &                          \\
\hline
\end{tabular}
\end{table}
}

\eat{
\begin{table}[h]
\small
\centering
\caption{Computational results of the ORLIB instances.}
\label{Tab:ORLIB}
\setlength{\tabcolsep}{1.0mm}{
\begin{tabular}{llllllllllllll}
\toprule
\multirow{2}{*}{Inst.} & \multicolumn{1}{c}{\multirow{2}{*}{$\delta$}} & \multicolumn{1}{c}{\multirow{2}{*}{P}} & \multicolumn{3}{c}{\fontfamily{qcr}\selectfont{SRM}}                                                          &      ~~                & \multicolumn{3}{c}{\fontfamily{qcr}\selectfont{Benders}}                                                       &      ~~                & \multicolumn{3}{c}{\fontfamily{qcr}\selectfont{Benders-AS}}                                                      \\
\cline{4-6}\cline{8-10}\cline{12-14}
                       & \multicolumn{1}{c}{}                          & \multicolumn{1}{c}{}                   & \multicolumn{1}{c}{t[s]} & \multicolumn{1}{c}{rgap[\%]} & \multicolumn{1}{c}{Obj} &                      & \multicolumn{1}{c}{t[s]} & \multicolumn{1}{c}{rgap[\%]} & \multicolumn{1}{c}{Obj} &                      & \multicolumn{1}{c}{t[s]} & \multicolumn{1}{c}{rgap[\%]} & \multicolumn{1}{c}{Obj}  \\

\hline
capa                   & 0.1                                           & 5                                      & 1858.3                   & 0                            & 10091227                &                                    & 1298.7                   & 0                            & 10091227                &                                    & \textbf{20.7 }           & 0                            & 10091227                 \\
                       & 0.1                                           & 10                                     & 2537.0                   & 0                            & 7082643                 &                                    & 1813.8                   & 0                            & 7082643                 &                                    & \textbf{26.7 }           & 0                            & 7082643                  \\
                       & 0.3                                           & 5                                      & 2400.9                   & 0                            & 10252897                &                                    & 2226.5                   & 0                            & 10252897                &                                    & \textbf{36.2 }           & 0                            & 10252897                 \\
                       & 0.3                                           & 10                                     & 3088.6                   & 0                            & 7175161                 &                                    & 2516.4                   & 0                            & 7175161                 &                                    & \textbf{72.5 }           & 0                            & 7175161                  \\
                       & 0.5                                           & 5                                      & 7200.0                   & 2.64                         & 10665607                &                                    & 4819.1                   & 0                            & 10665607                &                                    & \textbf{449.6 }          & 0                            & 10665607                 \\
                       & 0.5                                           & 10                                     & 7200.0                   & 4.27                         & 7502149                 &                                    & 7200.0                   & 2.24                         & 7530690                 &                                    & 7200.0                   & \textbf{0.63 }               & 7502149                  \\
                       &                                               &                                        &                          &                              &                         &                                    &                          &                              &                         &                                    &                          &                              &                          \\
capb                   & 0.1                                           & 5                                      & 2120.1                   & 0                            & 10206163                &                                    & 986.4                    & 0                            & 10206163                &                                    & \textbf{22.8 }           & 0                            & 10206163                 \\
                       & 0.1                                           & 10                                     & 1873.0                   & 0                            & 7037432                 &                                    & 1458.0                   & 0                            & 7037432                 &                                    & \textbf{22.3 }           & 0                            & 7037432                  \\
                       & 0.3                                           & 5                                      & 2740.5                   & 0                            & 10324019                &                                    & 1387.9                   & 0                            & 10324019                &                                    & \textbf{47.7 }           & 0                            & 10324019                 \\
                       & 0.3                                           & 10                                     & 2598.2                   & 0                            & 7162656                 &                                    & 2556.6                   & 0                            & 7162656                 &                                    & \textbf{50.1 }           & 0                            & 7162656                  \\
                       & 0.5                                           & 5                                      & 7200.0                   & 1.89                         & 10817993                &                                    & 5030.1                   & 0                            & 10818027                &                                    & \textbf{620.5 }          & 0                            & 10818027                 \\
                       & 0.5                                           & 10                                     & 7200.0                   & 3.51                         & 7467029                 &                                    & 6816.5                   & 0                            & 7467031                 &                                    & \textbf{2464.6 }         & 0                            & 7467031                  \\
                       &                                               &                                        &                          &                              &                         &                                    &                          &                              &                         &                                    &                          &                              &                          \\
capc                   & 0.1                                           & 5                                      & 1770.8                   & 0                            & 10075539                &                                    & 1499.1                   & 0                            & 10075539                &                                    & \textbf{24.1 }           & 0                            & 10075539                 \\
                       & 0.1                                           & 10                                     & 2307.1                   & 0                            & 6944169                 &                                    & 1531.4                   & 0                            & 6944169                 &                                    & \textbf{27.1 }           & 0                            & 6944169                  \\
                       & 0.3                                           & 5                                      & 2444.6                   & 0                            & 10230758                &                                    & 2198.8                   & 0                            & 10230758                &                                    & \textbf{61.2 }           & 0                            & 10230758                 \\
                       & 0.3                                           & 10                                     & 2411.4                   & 0                            & 7018687                 &                                    & 1519.2                   & 0                            & 7018687                 &                                    & \textbf{59.7 }           & 0                            & 7018687                  \\
                       & 0.5                                           & 5                                      & 7200.0                   & 4.60                         & 10700077                &                                    & 7200.0                   & 0.47                         & 10712575                &                                    & \textbf{1120.2 }         & 0                            & 10700077                 \\
                       & 0.5                                           & 10                                     & 7200.0                   & 5.09                         & 7448787                 &                                    & 7200.0                   & 3.28                         & 7474738                 &                                    & 7200.0                   & \textbf{0.89 }               & 7451464                  \\
\bottomrule

\end{tabular}}
\end{table}
}

}

\clearpage

\bibliographystyle{apalike}
\bibliography{ref}

\end{document}